\documentclass[10pt]{article}
\usepackage{amssymb}
\usepackage{amsmath}
 \usepackage{a4wide}
 \usepackage{cite}
 \usepackage{subcaption}

\usepackage{kbordermatrix}
 \usepackage{hhline}
 \usepackage{graphicx,latexsym,amscd,euscript,epsfig}
\newtheorem{theorem}{\bf Theorem}[section]
\newtheorem{proposition}[theorem]{\bf Proposition}
\newtheorem{definition}[theorem]{\bf Definition}
\newtheorem{corollary}[theorem]{\bf Corollary}

\newtheorem{exam}[theorem]{\bf Example}
\newtheorem{remark}[theorem]{\bf Remark}
\newtheorem{lemma}[theorem]{\bf Lemma}

\newcommand{\vone}{\vskip 2ex}

\newcommand{\dm}[1]{ {\displaystyle{#1} } }
\newcommand{\be}{\begin{equation}}
\newcommand{\ee}{\end{equation}}
\def\bmatrix#1{\left[ \begin{matrix} #1 \end{matrix} \right]}
\def \noin{\noindent}
 
\def \R{{\mathbb R}}
\def \C{{\mathbb C}}
\def \H{{\mathbb H}}

\def \lam{\lambda}

\def \Lam{\Lambda}

\def \diag{\mathrm{diag}}

\newcommand{\beano}{\begin{eqnarray*}}
\newcommand{\eeano}{\end{eqnarray*}}
\newcommand {\proof} {\par{\it Proof}. \ignorespaces}

\def \det{\mathrm{det}}

\renewcommand{\thefootnote}{\fnsymbol{footnote}}
\date{\today}

\title{Localization theorems for matrices and bounds for the zeros of  polynomials over a quaternion division algebra }
\author{Sk. Safique Ahmad\footnotemark[1] \and Istkhar Ali\footnotemark[2]}
\begin{document}
\maketitle

\begin{abstract} 
\noindent In this paper, Ostrowski and Brauer type theorems are derived for the left and right eigenvalues of a
quaternionic matrix. Generalizations of Gerschgorin type theorems are discussed for the
left and  the right eigenvalues of a quaternionic matrix. Thereafter a sufficient condition for the stability of a
quaternionic matrix is given that generalizes the stability condition for a complex matrix. Finally,
a characterization of bounds for the zeros of quaternionic polynomials is presented.
\end{abstract}
\vone
 \vone
\renewcommand{\thefootnote}{\fnsymbol{footnote}}
\footnotetext[1]{Corresponding author:  School of Basic Sciences, Discipline  of Mathematics, Indian Institute of Technology Indore, Simrol, Indore-453552, Madhya Pradesh, India.;
 \texttt{email: safique@iiti.ac.in}, Phone: +91-731-2438947, Fax: +91-731-2438933 }

 \footnotetext[2]{
 School of Basic Sciences, Discipline  of Mathematics, Indian Institute of Technology Indore, Simrol, Indore-453552, Madhya Pradesh, India, \texttt{email: istkhara@iiti.ac.in}.
 Research work funded by the CSIR, Govt. of India. }
\section{Introduction}\label{s1}
 This paper attempts to study localization theorems for matrices over a
 quaternion division algebra, which include the Ostrowski, Brauer, and 
 Gerschgorin type of theorems. Bounds for the zeros of quaternionic polynomials 
 are also considered. Localization theorems for quaternionic matrices have
 received much attention in the literature due to their applications in pure and applied 
 sciences, especially in quantum theory \cite{sla95, a99, jd02, arv89, t80, lw01, 
 gp02, l49, l14, rppv05, lr12, jw08, fz97, fz07, m99}. Unlike the case of matrices over the 
 field of complex numbers \cite{a46, s31, rc96, am37, r04}, localization theorems for
 quaternionic matrices have been proposed for left and right eigenvalues separately 
 in \cite{fz07,wzcl08, lyj12}. Ostrowski  and Brauer type theorems for the right 
 eigenvalues of a quaternionic matrix with all real diagonal entries have been introduced 
 in \cite{lyj12}. A Brauer type theorem for the left eigenvalues of a quaternionic matrix has
been considered in \cite[Theorem 4]{wzcl08}. Moreover, localization theorems for special 
quaternionic matrices, for instance, central closed quaternionic matrices, have been 
presented in \cite{wzcl08}.

In the first part of this paper, we provide a general framework for localization 
theorems for quaternionic matrices. Let $M_n(\H)$ be the space of all $n \times n$ 
quaternionic matrices. Then, for any $A=(a_{ij}) \in M_n(\H),$ we prove a Ostrowski type 
theorem which states that all the left eigenvalues of $A$ are located in the 
union of $n$ balls $T_i(A):= \{z \in \H:|z-a_{ii}| \leq r_i(A)^{\gamma} c_i(A)^{1-\gamma} \}$, where
$r_i(A):= \sum_{j=1,\, j\neq i}^n|a_{ij}|$ and $c_i(A):=
\sum_{j=1,\, j\neq i}^n |a_{ji}|,\,\,\forall\,\, \gamma \in [0,\,1]$. From this result, we deduce a
sufficient condition for invertibility of a quaternionic matrix. We also proved that the Ostrowski 
type theorem is valid for the right eigenvalues when all the diagonal entries of
the quaternionic matrix $A$ are real.

We find that the Brauer type theorem, proved in \cite[Theorem 5]{wzcl08} for the left eigenvalues in the case of deleted absolute column sums
of a quaternionic matrix, is incorrect, and we prove a corrected version. In addition, we derive some stronger results than \cite[Theorems 6, 7]{wzcl08} and \cite[Theorem 4.3]{lyj12}.
In fact, in the case of the generalized H$\ddot{\mbox{o}}$lder inequality over the skew field of quaternions, we show that all
the left eigenvalues of $A=(a_{ij}) \in M_n(\H)$ are contained in
the union of $n$ generalized balls: $B_i(A) :=\{z\in \H : |z-a_{ii}|\le (n-1)^{\frac{1- \gamma}{q}} r_i(A)^{\gamma}
(n_i^{(p)}(A))^{1- \gamma}\}$, where $\gamma \in [0,\, 1],$ $ n_i^{(p)}(A):= \left( \sum_{j=1,\, j\neq i}^{n} |a_{ij}|^{p}\right)^{\frac{1}
{p}}$, for any $ p, q \in (1,\,\infty)$ with $\frac{1}{p}+ \frac{1}{q} = 1$. Further,
we prove that all the right eigenvalues of $A\in M_n(\H)$ with all real diagonal entries are contained in the
union of $n$ generalized balls $B_i(A).$ In the sequel, we present localization theorems for the right
eigenvalues of quaternionic matrices.

In the second part of this paper, we provide bounds for the zeros of quaternionic polynomials 
using the aforementioned localization theorems. Recall that quaternionic polynomials in general 
are expressed in the following forms
\begin{eqnarray}
  p_l(z) &:=& q_m z^m+ q_{m-1} z^{m-1}+ \dots+q_1z+ q_0,\label{l}\\
 p_r(z) &:=& z^m q_m+ z^{m-1}q_{m-1}+ \dots+zq_1+ q_0,\label{lk}
\end{eqnarray}
where $q_{j},\,\,z \in \H,\,\, (0 \leq j \leq m).$ The polynomials ({\ref{l}}) and (\ref{lk}) are called simple
and monic if $q_m=1.$ Some recent developments on the location and computation of zeros of quaternionic polynomials
can be found in \cite{bt65, ddg10, dg10, sgv06, in41, go09, am04, rp01}. 

As a consequence of the localization theorems for quaternionic matrices, we provide sharper bounds compared to the bound
introduced by G. Opfer in \cite{go09} for the zeros of quaternionic polynomials. Finally, we provide bounds for the zeros of
quaternionic polynomials in terms of powers of the companion matrices associated with the quaternionic polynomials (\ref{l}) and (\ref{lk}). 
 Some of our bounds are sharper than the bound from \cite{go09}.

The paper is organized as follows: Section \ref{s2} reviews some existing results from \cite{fz97, rpp08}.
Section \ref{ss3} discusses the Greshgorin type, Ostrowski type, and Brauer type theorems for the left
and right eigenvalues of a quaternionic matrix. Section \ref{s4} explains bounds for the zeros of $p_l (z)$ and $p_r (z)$.
Comparisons are made with the bound provided in \cite{go09}. A sufficient condition for the stability of a
quaternionic matrix is also given. Section \ref{s6} introduces bounds for the zeros of the polynomials
$p_l (z)$ and $p_r (z)$ in terms of powers of their companion matrices.
Finally, Section \ref{cf} summarizes this work.
\section{Preliminaries}\label{s2}
{\bf Notation:} Throughout the paper, $\R$ and $\C$ denote the fields of real and complex numbers, respectively.
The set of real  quaternions is defined by  $$\H:= \left\{ q= a_0 + a_1{\bf{i}} +
a_2{\bf{j}} + a_3{\bf{k}}: a_0, a_1, a_2, a_3 \in \R \right\}$$ with ${\bf{i}}^2={\bf{j}}^2={\bf{k}}^2={\bf{ijk}}=-1.$
 The conjugate of $q \in \H$ is $\overline{q}:= a_0 - a_1{\bf{i}} - a_2{\bf{j}} - a_3{\bf{k}}$  and the modulus of $q$ is
$|q|: =  \sqrt {a_0^2 + a_1^2 + a_2^2 + a_3^2}$.
$\Im{(a)}$ denotes the imaginary part of $a \in \C$. The real part of a quaternion $q= a_0 + a_1{\bf{i}} +
a_2{\bf{j}} + a_3{\bf{k}}$ is defined as 
$\Re(q)=a_0.$
The collection of all $n$-column vectors with elements in $\H$ is denoted by $\H^{n}$.
For $x \in \mathcal{K}^n,$ where $\mathcal{K} \in \{ \R, \C, \H\},$ the transpose of $x$ is $x^T.$ If $x=[x_1,\ldots,x_n]^T,$ the conjugate of
$x$ is defined as $\overline{x}=[\overline{x_1}, \ldots, \overline{x_n}]^T$ and the conjugate transpose of $x$
is defined as $x^H=[\overline{x_1}, \ldots, \overline{x_n}].$
For $x, y \in \H^n,$ the inner product is defined as $\langle x, y \rangle:= y^Hx$ and the norm of $x$ is defined as
$\|x\|:= \sqrt{\langle x, x\rangle}$.
The sets of $m\times n$  real, complex, and quaternionic matrices are denoted by $M_{m \times n}(\R),$ $ M_{m \times n}(\C),$
and $M_{m \times n}(\H),$ respectively. When $m=n$, these sets are
denoted by $M_n(\mathcal{K})$, $\mathcal{K} \in \{\R, \C, \H \}$.
%
For $A \in M_{m \times n}({\mathcal{K}}),$ 
the conjugate, transpose, and conjugate transpose 
 of $A$ are defined as $\overline{A}=(\overline{a_{ij}})$, $A^T=(a_{ji}) \in M_{n \times m} (\H),
 $ and $A^H=(\overline{A})^T \in M_{n \times m}(\H),$ respectively.
 For $z\in \H^n,$ the vector $p$-norm on $\H^n$ is
defined by  $\|z\|_p:=(\sum_{i=1}^n|z_i|^p)^{1/p},$ where $1\le p < \infty$ and $\|z\|_\infty:=\dm {\max_{1\le i \le n}}\{|z_i|\}.$ Define
$\R^{+}:=\{ \alpha: \alpha \in \R, \alpha > 0 \}.$
The set 
 \begin{eqnarray*}
  [q]:=\{ r \in \H : r=\rho^{-1} \, q \,\rho \,\, for\,\, all\,\, 0\not=\rho \in \H\}
 \end{eqnarray*}
 is called an equivalence class of $q \in \H.$
%
Let $x\in \H^n$. Then $x$ can be uniquely expressed as $x= x_1+ x_2{\bf{j}},$  where $x_1, x_2 \in \C^n.$ Define the function $ \psi :  \H^n \rightarrow
\C^{2n}$ by 
$$ \psi_x:=\bmatrix{x_1 \\-\overline{x_2}}.$$
 This function $\psi$ is an injective linear transformation from $\H^n$ to $\C^{2n}.$
\begin{definition}\label{ch1cam}
Let $A\in M_n({\H})$. Then $A$ can be uniquely expressed as $A= A_1+ A_2{\bf{j}},$  where $A_1,A_2 \in M_n(\C).$ Define the function
$\Psi :  M_n({\H}) \rightarrow
M_{2n}({\C})$ by 
\begin{eqnarray*}
\Psi_A:=\bmatrix{A_1 & A_2\\-\overline{A_2}&\overline{A_1}}.
\end{eqnarray*}
The matrix $\Psi_A$ is called 
the complex adjoint matrix of $A$.
\end{definition}
%
\begin{definition}\label{def3.1}
Let  $A\in M_n({\H}).$ Then the left, right, and the standard eigenvalues, respectively, are given by
\begin{eqnarray*}
\Lambda_l(A) &:=& \left\{\lam \in \H : Ax = \lam x \,\, \mbox{for some nonzero}\,\,
x \in \H^{n} \right\},\\
\Lambda_r(A) &:=& \left\{\lam \in \H : Ax = x \lam \,\, \mbox{for some nonzero}\,
\, x \in \H^{n} \right\}\,\mbox{and}\\
\Lambda_s(A) &:=& \left\{\lam \in \C : Ax =  x\lam \,\, \mbox{for some nonzero}\,
 \, x \in \H^{n},\, \Im(\lam ) \geq 0 \right\}.
 \end{eqnarray*}
 \end{definition}
\begin{definition}
 Let $A \in M_n(\H).$ Then $A$ is said to be a central closed matrix if there exists an invertible matrix $T$ such that
 \[
  T^{-1} A T=\diag(\lam_1,\lam_2,\ldots,\lam_n),\quad
\mbox{where}\,\, \lam_i \in \R, \quad 1 \leq i \leq n. \]
\end{definition}
\begin{definition}\label{st}
   Let $A\in M_n(\H).$ Then the matrix $A$ is said to be stable if and only if $\Lam_r(A) \subset \H^{-}:= \left\{ q \in \H :
   \Re(q) < 0\right\}.$
  \end{definition}
 \begin{definition}
 Let $A \in M_n(\H).$ Then $A$ is said to be $\eta$-Hermitian if $A= \left(A^{\eta}\right)^H,$
 where $A^{\eta}=\eta^H A \eta $ and  $\eta \in \{{\bf{i}}, {\bf{j}}, {\bf{k}} \}.$
 \end{definition}
 \begin{definition}
 A matrix $A \in M_n(\H)$ is said to be invertible if there exists
 $B \in M_n(\H)$ such that $AB=BA=I_n,$ where $I_n$ is the $n \times n$ identity matrix.
 \end{definition}
We next recall the following result necessary for the development of our theory.
%
 %
%
\begin{theorem}\cite[Theorem 4.3]{fz97}.\label{nst}
 Let $A\in M_n({\H}).$ Then the following statements are equivalent:

\noin $(a)$ $A$ is invertible,\,\, $(b)$ $Ax = 0$ has the unique solution, \,\,
 $(c)$  $\det(\Psi_A) \neq 0,$\,\,
 $(d)$  $\Psi_A$ is invertible,\,\,
 $(e)$ $A$ has no zero eigenvalue.
\end{theorem}
%

Let $A:= (a_{ij})\in M_n(\H)$ and define the absolute row and column sums of $A$ as
 \[
 r_i'(A):=r_i(A)+ |a_{ii}|\,\,\mbox{and}\,\, c_i'(A):=c_i(A)+|a_{ii}|  \,\,\,\,(1 \leq i \leq n).
 \]
\section{Distribution of the left and right eigenvalues of quaternionic matrices}\label{ss3}

It is known from \cite[Corollary 3.2]{l08} that a quaternionic matrix $A$ and its conjugate transpose $A^H$ have the same right
eigenvalues. However, $A$ and $A^H$ may not have the same left eigenvalues, take for example
                        $A= \bmatrix{{\bf{i}}& 0 \\0 & {\bf{j}}}$ and $A^H= \bmatrix{{\bf{-i}}
                        & 0 \\0 &{\bf{-j}}}$. We now present the following lemma for left eigenvalues of $A$ and $A^H.$
\begin{lemma}\label{prop4}
Let $A\in M_n({\H})$ and let $\lam \in \H.$ Then $\lam$ is a left eigenvalue of $A$ if and only if $\overline{\lam}$
is a left eigenvalue of  $A^H.$
\end{lemma}
\proof Let $\lam$ be a left eigenvalue of $A.$ Then there exists
$ x (  \neq 0) \in \H^n$ such that  $(A-\lam I_n)x= 0.$ This can be written as $\Psi_{(A-\lam I_n)}
\psi_x=0.$ Hence it follows that $\lam$ is a left eigenvalue of $A$ if and only if
$\det\left[\Psi_{(A-\lam I_n)}\right]=0$ $\Leftrightarrow \det\left[\Psi^H_{(A-\lam I_n)}\right]=0
\Leftrightarrow  \det\left[\Psi_{(A-\lam I_n)^H}\right]=0
\Leftrightarrow  \det\left[\Psi_{(A^H-\overline{\lam} I_n)}\right]=0.
$
Thus, $\overline{\lam}$ is a left eigenvalue of $A^H.\,\,\,\blacksquare$

 The Gerschgorin type theorem for the left eigenvalues using deleted absolute row sums of a matrix $A\in M_n (\H)$ is proved in \cite{fz07}.
 However, the Gerschgorin type theorem
for the left eigenvalues using deleted absolute column sums of $A$ has not yet been established.
We now state and prove the theorem.
\begin{theorem}\label{thm2}
 Let $A:=(a_{ij}) \in M_n({\H}).$ Then all the left eigenvalues of $A$ are
 located  in the union of $n$ Gerschgorin balls $\Omega_i(A):=\left\{ z \in \H:
 |z-a_{ii}|\leq c_i(A) \right\},  1 \leq i \leq n,$ that is,
 \[
\Lam_l(A) \subseteq \Omega(A):=\cup_{i=1}^{n} \Omega_i(A).
\]
\end{theorem}
\proof  Let $\lam$ be a left eigenvalue of $A.$ Then from Lemma \ref{prop4},
$\overline{\lam}$ is a left eigenvalue of $A^H.$ Then there exists some nonzero $x \in \H^n$ such that
$A^H x = \overline{\lam}x$. Let $x:= [x_1,\ldots, x_n]^T \in \H^n$ and let
$x_t$ be an element of $x$ such that $|x_t|\ge |x_i|, 1\le i \le n$.
Then, $|x_t| > 0$. From the $t$-th equation of $A^Hx= \overline{\lam} x$,  we have
\begin{eqnarray*}
 \sum_{j=1}^n \overline{a_{jt}} x_j &=& \overline{\lam} x_t.
\end{eqnarray*}
This shows
\begin{eqnarray*}
|\lam- a_{tt}| &\le & \sum_{j=1,\, j\neq t}^n |a_{jt}|:=c_t(A).\,\,\,\blacksquare
\end{eqnarray*}

We now have the following localization theorem for the deleted absolute
 row and column sums of a matrix $A \in M_n(\H)$ which is known as {\em Ostrowski type
 theorem}.
%
\begin{theorem}$($Ostrowski type theorem for the left eigenvalues$)$ \label{os1}
Let $A:=(a_{ij})\in M_n({\H})$ and let $\gamma \in [0, 1]$. Then all the
left eigenvalues of $A$ are located in the  union of $n$ balls $T_i(A):=
\{z \in \H:|z-a_{ii}| \leq r_i(A)^{\gamma} c_i(A)^{1-\gamma}\}, 1 \leq i \leq n,$  that is,
 \[
  \Lam_l(A) \subseteq T(A):= \cup_{i=1}^{n} T_i(A).
  \]
\end{theorem}
\proof Let $\lam$ be a left eigenvalue of $A.$ Then by \cite[Theorem $6$]{fz07}, for $\gamma\in [0,\,1],$ we have
\begin{equation}\label{eqbp}
 |\lambda-a_{ii}|^{\gamma} \le r_i(A)^{\gamma}, \quad 1 \leq i \leq n.
\end{equation}
Similarly, from Theorem \ref{thm2}, we obtain
\begin{equation}\label{eqbq}
 |\lambda-a_{ii}|^{1-\gamma}  \le c_i(A)^{1-\gamma}, \quad 1 \leq i \leq n.
\end{equation}
Combining (\ref{eqbp}) and (\ref{eqbq}), we get
\[
 |\lambda- a_{ii}| \le r_i(A)^{\gamma}  c_i(A)^{1-\gamma}, \quad 1 \leq i \leq n.
\]
Thus, all the left eigenvalues of $A$ are located in the union of $n$ balls $T_i(A)$. $\blacksquare$

Next, we derive Ostrowski type theorem for right eigenvalues of $A\in M_n(\H)$ with all real diagonal entries.

\begin{theorem}\label{j12}
Let $A:=(a_{ij}) \in M_n({\H})$ with $a_{ii}\in \R$ and let $\gamma \in [0,\,1]$. Then all the right eigenvalues of $A$ are located in the union of $n$
balls $G_i(A):=\left\{z\in \H : |z-a_{ii}|\le r_i(A)^{\gamma}
c_i(A)^{1-\gamma}\right\}$, $1 \leq i \leq n,$  that is,
\[ \Lambda_r(A)\subseteq G(A):= \cup_{i=1}^{n} G_i(A).\]
\end{theorem}
\proof Let $\lam$ be a right eigenvalue of $A.$ Then there exists some nonzero $x \in \H^n$ such that
$A x = x \lam $. Let $x:= [x_1,\ldots, x_n]^T \in \H^n$ and let
$x_t$ be an element of $x$ such that $|x_t|\ge |x_i|, 1\le i \le n$.
From the $t$-th equation of $Ax= x \lam $,  we have
\begin{eqnarray}
 a_{tt} x_t +\sum_{j=1,\, j\neq t}^n a_{tj} x_j &=& x_t \lam.
\end{eqnarray}
Since $a_{tt} \in \R,$ $a_{tt} x_t=x_t a_{tt}.$ Proceeding as in the proof of Theorem \ref{thm2}, we obtain
    \begin{eqnarray}\label{R1}
     |\lam- a_{tt}| \le  \sum_{j=1,\, j\neq t}^n |a_{jt}|=:r_t(A).
    \end{eqnarray}
From \cite[Corollary 2.7]{l08}, $\lam$ is also a right eigenvalue of $A^H$. Then
 \begin{eqnarray}\label{R2}
     |\lam- a_{tt}| \le  \sum_{j=1,\, j\neq t}^n |a_{tj}|=:c_t(A).
\end{eqnarray}
Let $\gamma \in [0, 1].$  Then from (\ref{R1}) and   (\ref{R2}), we obtain
 \begin{eqnarray}\label{R3}
          |\lam- a_{tt}|^{\gamma} \le  r^{\gamma}_t(A),
          \end{eqnarray}
   \begin{eqnarray}\label{R4}
 |\lam- a_{tt}|^{1-\gamma} \le  c^{1-\gamma}_t(A).
 \end{eqnarray}
  Combining (\ref{R3}) and (\ref{R4}), we get
  \[
    |\lam- a_{tt}| \le r(A)^{\gamma}_t \, c(A)^{1-\gamma}_t.\,\,\,\blacksquare
  \]

\begin{corollary}\label{s3}
 For any  $A:= (a_{ij})\in  M_n(\H),$ $n\geq 2$ and for any $ \gamma \in [0,\,1]$. Let us assume that
\begin{equation}\label{eqn5}
|a_{ii}| > r_i(A)^{\gamma} \ r_i(A)^{1-\gamma}, \quad 1\le i\le n.
\end{equation}
 Then $A$ is invertible.
 \end{corollary}
 \proof On the contrary, suppose $A$ is not invertible. Then by Theorem \ref{nst},
 there is a left eigenvalue $\lam=0$ of  $A$. Now from Theorem \ref{os1}, we obtain
 $|a_{ii}| \leq r_i(A)^{\gamma}  c_i(A)^{1-\gamma}$. This contradicts our assumption
 (\ref{eqn5}). Hence $A$ is invertible.$\,\,\,\blacksquare$
%
%
%

It is known that a quaternionic matrix $A \in M_n(\H)$ may have at most $2n$ complex right eigenvalues. From
Theorem \ref{j12}, all the complex right eigenvalues of a matrix $A=(a_{ij}) \in M_n(\H)$
with all real diagonal entries lie in the union of $n$-discs $\mathcal{E}_i(A):=\{z \in \C: |z-a_{ii}|
\leq r_i(A)^{\gamma} c_i(A)^{1-\gamma}  \},$ $1 \leq i \leq n,$  that is,
    \begin{eqnarray}\label{EEE2}
   \Lambda_c(A)\subseteq \mathcal{E}(A):= \cup_{i= 1}^{n}\mathcal{E}_i(A),\,
   \mbox{where}\,\, \Lambda_c(A):=\{\lam \in \C : A x=x \lam,\,\, 0 \neq x \in \H^n \}.
 \end{eqnarray}
%

The Brauer type theorem is proved in \cite{wzcl08} for the left eigenvalues in the case of deleted absolute
column sums of a matrix $A \in M_n(\H).$ That is, if $\lam \in \Lam_l (A),$ then its conjugate $\overline{\lam}$ lies in the
union of $\frac{n (n- 1)}{2}$ ovals of Cassini. However, this is incorrect as the following example suggest:

\begin{exam}\label{count2}{\rm
 Let $A=\bmatrix{{\bf{i}} &{\bf{k}} \\ 0 & {\bf{j}}}.$
 Then by \cite[Theorem 5]{wzcl08}, oval of Cassini is given by
 $\left\{z \in \H : |z-{\bf{i}}|\ |z-{\bf{j}}|\leq 0\right\}.$ Here,
 ${\bf i}$ is a left eigenvalue of $A$ and its conjugate
 ${\bf -i}$  is not contained in the above oval of Cassini.}
 \end{exam}
 According to \cite[Theorem 5]{wzcl08}, if $\lam \in \Lam_l(A),$ then $\overline{\lam} \in \dm{\cup_{\substack{i, j=
1,\\ i \ne j}}^{n}} F_{ij}(A),$ where
$$F_{ij}(A):= \left\{z\in \H : |z-a_{ii}| \ |z-a_{jj}| \le c_i(A)
c_j(A)\right\},\quad 1 \leq i,j \leq n,\quad i \neq j.$$
However, this result is not necessarily true as
\[
 |\overline{\lam}-a_{ii}| \ |\overline{\lam}-a_{jj}| > c_i(A)
c_j(A),\quad 1 \leq i,j \leq n,\quad i \neq j,
\]
which follows from Example \ref{count2}. Now, we derive a corrected version of \cite[Theorem 5]{wzcl08} as follows:
\begin{theorem}\label{c}
 Let $A:=(a_{ij}) \in M_n({\H}).$ Then all the left
eigenvalues of $A$ are located in the union of $\frac{n(n-1)}{2}$ ovals of
Cassini
$$F_{ij}(A):= \left\{z\in \H : |z-a_{ii}| \ |z-a_{jj}| \le c_i(A) c_j(A) \right\},
\quad 1 \leq i,j \leq n, \quad i \neq j,$$
that is, $ \Lambda_l(A)\subseteq F(A):=\dm{\cup_{\substack{i, j=
1,\\ i \ne j}}^{n}} F_{ij}(A).$
\end{theorem}
\proof Let $\lam$ be a left eigenvalue of $A.$ Then by Lemma \ref{prop4},
$\overline{\lam}$ is a left eigenvalue of $A^H$. Then there exists some nonzero $x \in \H^n$ such that
$A^H x = \overline{\lam}x$. Let $x:= [x_1,\ldots, x_n]^T \in \H^n$ and let
$x_s$ be an element of $x$ such that 
$|x_s| \geq |x_i|$, $1 \leq i \leq n.$ Then, $|x_s| >0.$ Clearly, if all the other elements of $x$ are zero,
then the required result holds.

Let $x_s$ and $x_t$ be two nonzero elements of $x$ such
that $|x_s|\ge|x_t|\ge|x_i|, 1\le i \le n, i \neq s$.
From the $s$-th equation of $A^H x= \overline{\lam} x $, we have
 \[
 \sum_{j=1}^{n} \overline{ a_{js}} x_j= \overline{\lam} x_s,
 \]
which implies
 \[
(\overline{\lambda}-\overline{a_{ss}}) x_s= \sum_{j=1,\, j\neq s}^{n} \overline{a_{js}} x_j.
 \]
Thus
\begin{equation}\label{eqb1}
 |\lambda-a_{ss}| \le \left(\frac{|x_t|}{|x_s|}\right) \ c_s(A).
 \end{equation}
 Similarly, from $A^H x= \overline{\lam}x,$ we obtain
\begin{equation}\label{eqb2}
 |\lambda-a_{tt}| \le \left(\frac{|x_s|}{|x_t|}\right) \ c_t(A).
 \end{equation}
Combining (\ref{eqb1}) and (\ref{eqb2}), we have
 \[
|\lambda-a_{ss}| \ |\lambda-a_{tt}|\le c_s(A) c_t(A).
\]
Hence, all the left eigenvalues of $A$ are located in the union of $\frac{n(n-1)}{2}$
ovals of Cassini $F_{ij}(A), \quad 1 \leq i,j \leq n, \quad i \neq j.$\,\,\, $\blacksquare$

Theorem $7$ of \cite{wzcl08} was stated for a central closed quaternionic matrix. Now we
generalize this result for all quaternionic matrices as follows.

\begin{theorem}\label{rc}
Let $A:=(a_{ij}) \in M_n(\H)$ and let $\gamma \in [0,\,1]$. Then all the left
eigenvalues of $A$ are located in the union of $\dm{\frac{n(n-1)}{2}}$ ovals of Cassini
\[K_{ij}(A):=\left\{z\in \H : |z-a_{ii}| \ |z-a_{jj}|\le r_i(A)^{\gamma} \
r_j(A)^{\gamma} \ c_i(A)^{1-\gamma } \ c_j(A)^{1-\gamma} \right\}, \quad 1 \leq i,j \leq n, \quad i \neq j,\]
that is,
\[
 \Lambda_l(A)\subseteq K(A):=\dm{\cup_{\substack{i,j=1\\ i \ne j}}^{n}}K_{ij}(A).
 \]
 \end{theorem}
 \proof Let $\lam$ be a left eigenvalue of $A.$ Then by \cite[Theorem 4]{wzcl08} and Theorem \ref{c},
 for $\gamma \in [0, 1]$, we have
\begin{equation}\label{ch2eqn26}
 |\lambda-a_{ii}|^{\gamma} |\lambda-a_{jj}|^{\gamma} \le r_i(A)^{\gamma}
 r_j(A)^{\gamma}, \quad 1 \leq i,j \leq n, \quad i \neq j
\end{equation}
and
\begin{equation}\label{ch2eqn27}
 |\lambda-a_{ii}|^{1-\gamma} |\lambda-a_{jj}|^{1-\gamma}  \le c_i(A)^{1-\gamma} 
 c_j(A)^{1-\gamma}, \quad 1 \leq i,j \leq n, \quad i \neq j.
\end{equation}
Combining (\ref{ch2eqn26}) and (\ref{ch2eqn27}), we have
\[
 |\lambda- a_{ii}| |\lambda-a_{jj}| \le r_i(A)^{\gamma} r_j(A)^{\gamma} c_i(A)^{1-\gamma} c_j(A)^{1-\gamma},
  \quad 1 \leq i,j \leq n, \quad i \neq j.\,\,\, \blacksquare
\]
\begin{corollary}
For any  $A:= (a_{ij})\in  M_n(\H),$ $n\geq 2$ and for any $ \gamma \in [0,\,1]$. Assume that
%
$$|a_{ii}||a_{jj}|> r_i(A)^{\gamma}
r_j(A)^{\gamma} \ c_i(A)^{1-\gamma }  c_j(A)^{1-\gamma},\quad  1 \leq i,j \leq n,\quad i \neq j.$$
Then $A$ is invertible.
\end{corollary}
\begin{corollary}
 Let $A:= (a_{ij}) \in M_n({\H}).$ Then all the left
eigenvalues of $A$ are located in the union of $\frac{n(n-1)}{2}$ ovals of Cassini
$$ \Lambda_l(A)\subseteq \Phi(A):=\cup_{{\substack{i,j=1\\ i \ne j}}}^{n}
\left\{z\in \H : |z-a_{ii}| \ |z-a_{jj}| \le \min \{r_i(A) r_j(A),\, c_i(A) c_j(A) \} \right\}.$$
\end{corollary}
\proof Substituting $\gamma=0,1$ in Theorem \ref{rc}, we obtain the following:
\begin{itemize}
\item[$(a)$] $ \Lambda_l(A)\subseteq E(A):=\cup_{{\substack{i,j=1\\ i \ne j}}}^{n}
\left\{z\in \H : |z-a_{ii}| \ |z-a_{jj}| \le c_i(A)  c_j(A)\right\}.$
\item[$(b)$] $ \Lambda_l(A)\subseteq F(A):=\cup_{{\substack{i,j=1\\ i \ne j}}}^{n}
\left\{z\in \H : |z-a_{ii}| \ |z-a_{jj}| \le r_i(A) r_j(A)\right\}.$
\end{itemize}
Combining $(a)$ and $(b),$ we get the required result.$\,\,\,\blacksquare$

The following result provides better estimate than Theorem \ref{j12}.
\begin{theorem}\label{btt}
Let $A:= (a_{ij}) \in M_n(\H)$ with $a_{ii} \in \R$ and let $\gamma \in [0,1]$. Then all the right
eigenvalues of $A$ are located in the union of $\dm{\frac{n(n-1)}{2}}$ ovals of Cassini $\mathcal{G}_{ij}(A):=
\left\{z\in \H : |z-a_{ii}| \ |z-a_{jj}|
\le r_i(A)^{\gamma} \ r_j(A)^{\gamma} \ {c_i}(A)^{1-\gamma}\ {c_j}(A)^{1-\gamma}\right\}
, \quad 1 \leq i,j \leq n, \quad i \neq j,$ that is,
\[
\Lambda_r(A)\subseteq \mathcal{G}(A):=\cup_{\substack{i,j=1\\ i \ne j}}^{n}\mathcal{G}_{ij}(A).
\]
  \end{theorem}
\proof
Let $\lam$ be a right eigenvalue of $A.$ Then by \cite[Theorem 4.1, Corollary 4.1]{lyj12},
 for $\gamma \in [0, 1]$, we have
\begin{equation}\label{nch2eqn26}
 |\lambda-a_{ii}|^{\gamma} |\lambda-a_{jj}|^{\gamma} \le r_i(A)^{\gamma} 
 r_j(A)^{\gamma},\quad 1 \leq i,j \leq n, \quad i \neq j
\end{equation}
and
\begin{equation}\label{nch2eqn27}
 |\lambda-a_{ii}|^{1-\gamma} |\lambda-a_{jj}|^{1-\gamma}  \le c_i(A)^{1-\gamma}
 c_j(A)^{1-\gamma}, \quad 1 \leq i,j \leq n, \quad i \neq j.
\end{equation}
Combining (\ref{nch2eqn26}) and (\ref{nch2eqn27}), we have
\[
 |\lambda- a_{ii}| |\lambda-a_{jj}| \le r_i(A)^{\gamma} r_j(A)^{\gamma} c_i(A)^{1-\gamma} c_j(A)^{1-\gamma}
 ,\quad 1 \leq i,j \leq n, \quad i \neq j.\,\,\, \blacksquare
\]

From Theorem \ref{btt}, all the complex right eigenvalues of a matrix $A:=(a_{ij}) \in M_n(\H)$ with
$a_{ii} \in \R$, $1 \leq i \leq n$ are contained in the union of $\frac{n(n-1)}{n}$ ovals of Cassini
$\mathcal{F}_{ij}(A):=\{z \in \C: |z-a_{ii}|\, |z-a_{jj}| \leq r_i(A)^{\gamma} r_j(A)^{\gamma} \ {c_i}(A)^{1-\gamma}
{c_j}(A)^{1-\gamma} \}$, $ \quad 1\le i, j\le n, \quad i \ne j$,
 that is,
    \begin{eqnarray}\label{EEE3}
   \Lambda_c(A)\subseteq \mathcal{F}(A):= \cup_{\substack{i,j=1\\ i \ne j}}^{n} \mathcal{F}_{ij}(A),
   \end{eqnarray}

The following theorem shows that Theorem \ref{rc} is sharper than Theorem \ref{os1}.
\begin{theorem}\label{sh}
Let $A:= (a_{ij}) \in M_n(\H)$ with $ n \geq 2$ and let $\gamma \in [0, 1]$. Then
$$K(A)\subseteq T(A),$$
 where $G(A)$ and $\mathcal{G}(A)$ are defined in Theorem \ref{os1}  and Theorem \ref{rc}, respectively.
 \end{theorem}
\proof Let $z \in K_{ij}(A)$ and fix any $i$ and $j,\, \,( 1\le i,j \le n,\, i \neq j)$. Then from
 Theorem \ref{rc}, we have
\begin{equation}\label{Eqn23}
 |z-a_{ii}| \ |z-a_{jj}| \le r_i(A)^{\dm{\gamma}} r_j(A)^{\dm{\gamma}} c_i(A)^{\dm{1-\gamma}} c_j(A)^{\dm{1-\gamma}}.
 \end{equation}
 Now the following two cases are possible.

{\bf Case 1:} If $r_i(A)^{\dm{\gamma}} \  r_j(A)^{\dm{\gamma}} c_i(A)^{\dm{1-\gamma}} c_j(A)^{\dm{1-\gamma}}=0,$
then $z= a_{ii}$ or $z= a_{jj}$.
However, from Theorem \ref{os1}, we have $a_{ii}\in T_i(A)$
and $a_{jj} \in T_j(A).$ Thus $z \in T_i(A) \cup T_j(A).$

{\bf Case 2:}
If $r_i(A)^{\dm{\gamma}}  r_j(A)^{\dm{\gamma}}  c_i(A)^{\dm{1-\gamma}} c_j(A)^{\dm{1-\gamma}} > 0,$ then by (\ref{Eqn23})
\begin{equation}\label{Eqn24}
\left(\frac{|z-a_{ii}|}{r_i(A)^{\dm{\gamma}} c_i(A)^{\dm{1-\gamma}}}\right) \left(\frac{|z-a_{jj}|}{r_j(A)^{\dm{\gamma}}
c_j(A)^{\dm{1-\gamma}}} \right) \le 1.
\end{equation}
 As the left side of (\ref{Eqn24}) cannot exceed unity, one of the 
 factors of the left side can be at most unity, that is, $z \in T_i(A)$ or 
 $z \in T_j(A).$ Hence $z \in T_i(A) \cup  T_j(A)$. Thus
 \begin{eqnarray}\label{eqn25}
 K_{ij}\subseteq T_i(A) \cup T_j(A).
 \end{eqnarray}
 From Theorem \ref{os1} and Theorem \ref{rc}, we obtain
 \[K(A):= \cup_{\substack{i, j= 1\\ i\ne j}}^{n} K_{ij}(A) \subseteq \cup_{\substack{i, j= 1\\ i\ne j}}^{n} \left\{T_i(A) \cup T_j(A) \right\}= \cup_{k= 1}^n T_k(A)=: T(A).\,\,\, \blacksquare\]

Similarly, we have the following relation between Theorem \ref{btt} and Theorem \ref{j12}.
\begin{theorem}\label{shh}
Let $A:= (a_{ij})\in M_n(\H), n\ge 2$ with $ a_{ii}\in \R$ and let $\gamma\in [0, 1].$ Then
$$\mathcal{G}(A)\subseteq G(A),$$ where $G(A)$ and $\mathcal{G}(A)$ are defined in Theorem \ref{j12}
and Theorem \ref{btt}, respectively.
\end{theorem}
\proof The proof is similar to the proof of Theorem \ref{sh}.$\,\,\,\blacksquare$


The following example illustrates Theorem \ref{shh} for complex right eigenvalues of a matrix $A:=(a_{ij}) \in M_n(\H)$
with $a_{ii} \in \R$, $1 \leq i \leq n.$
\begin{exam}\label{ree1}{\em Let
$A= \bmatrix{3 & 1+{\bf{i}}+{\bf{j}}-{\bf{k}} & 2+3{\bf{j}}-\sqrt{3} {\bf{k}} \\
5+ \sqrt{2}{\bf{j}}+3 {\bf{k}} & -2 & 3{\bf{j}}+4{\bf{k}} \\
4+3{\bf{j}} & 2-{\bf{i}}-2{\bf{k}} & -5 }.$
Substituting $\gamma=1/4$ in (\ref{EEE2}), we get the following three discs:
%
\begin{eqnarray*}
 \mathcal{E}_1(A)&:=&\{z \in \C: |z-3| \leq 9.4533 \},\\ \mathcal{E}_2(A)&:=&\{z \in \C: |z+2| \leq 6.0894\},\\
 \mathcal{E}_3(A)&:=&\{z \in \C: |z+5| \leq 8.7389\}.
\end{eqnarray*}
Similarly, let $\gamma=1/4$ in (\ref{EEE3}), we get the following three discs:
\begin{eqnarray*}
 \mathcal{F}_{12}(A) &:=& \{z \in \C : |z-3|\, |z+2| \leq 57.5649 \},\\
 \mathcal{F}_{23}(A) &:=& \{z \in \C : |z+2|\, |z+5| \leq 53.2145 \},\\
 \mathcal{F}_{31}(A) &:=& \{z \in \C : |z+5|\, |z-3| \leq 82.6108 \}.
 \end{eqnarray*}
In this example, there are six complex right eigenvalues $\lam_j\,\, ( 1 \leq  j \leq  6)$ 
which are shown in Figure \ref{fig1}.
The set $\mathcal{F}(A):= \mathcal{F}_{12}(A) \cup \mathcal{F}_{23}(A) \cup \mathcal{F}_{31}(A)$ is
represented by shaded region in Figure \ref{fig1}. From Figure \ref{fig1}, it is clear that
$\mathcal{F}(A) \subset \mathcal{E}(A),$ where $\mathcal{E}(A):=\mathcal{E}_1(A) \cup \mathcal{E}_2(A)
\cup \mathcal{E}_3(A).$
 \begin{figure}[h!]
  \hspace{-1.3cm}\includegraphics[height=8cm,width=15cm]{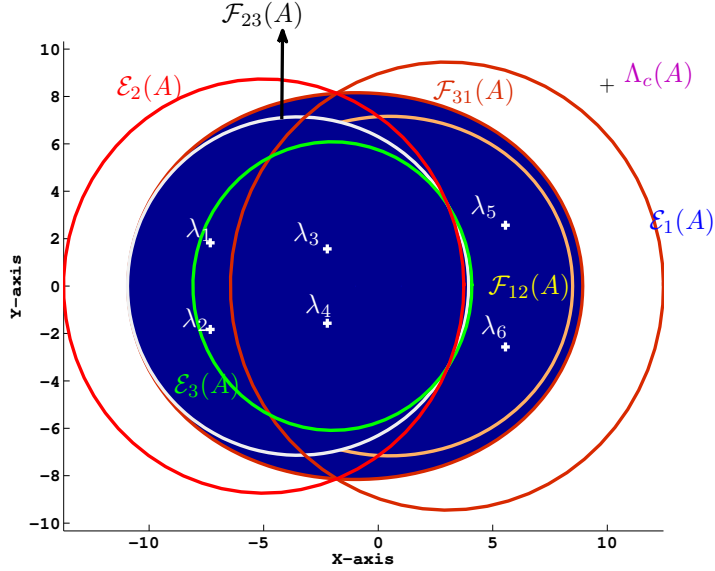}
   \caption{Location of the complex right eigenvalues of the matrix $A$ from Example \ref{ree1}.}\label{fig1}
  \end{figure}
}
 \end{exam}
%
 For $A:= (a_{ij})\in M_n(\H)$, define
 $$n_i^{(p)}(A):= \left( \sum_{j=1,\, j\neq i}^{n} |a_{ij}|^{p}\right)^{\frac{1}{p}}, \quad
 1 \leq i \leq n, \quad p \in (1,\infty).$$

We are now ready to derive the following localization theorem for left eigenvalues of a quaternionic matrix.
 \begin{theorem}\label{g1}
  Let $A:=(a_{ij}) \in M_n({\H})$ and let $\gamma \in [0, 1]$. Then all the left eigenvalues of $A$ are  contained in the
  union of  $n$ generalized balls \[B_i(A):= \left\{z\in \H : |z- a_{ii}| \le
 (n-1)^{\frac{1-\gamma}{q}} r_i(A)^{\gamma} (n_i^{(p)}(A))^{1- \gamma}\right\},\quad 1 \leq i \leq n,\]
 that is,
 \[\Lam_l(A) \subseteq B(A):= \cup_{i=1}^{n}B_i(A),\]
 for any $p, q \in (1, \infty)$ with $\frac{1}{p}+\frac{1}{q}=1.$
 \end{theorem}
  \proof Let $\mu$ be a left eigenvalue of $A.$ Then there exists some nonzero $x \in \H^n$ such that
$A x = {\mu}x$. Let $x:= [x_1,\ldots, x_n]^T \in \H^n$ and let
$x_t$ be an element of $x$ such that $|x_t|\ge |x_i|,$ $1\le i \le n$.
Then from $Ax= \mu x,$ we have
\begin{eqnarray*}
 a_{tt} x_t + \dm{\sum_{ j= 1,\, j \neq t}^{n} a_{tj} x_j} &=&  \mu x_t.
\end{eqnarray*}
This implies
\begin{equation}\label{eqn9}
 |\mu-a_{tt}||x_t|= \left|\sum_{j= 1,\, j \neq t}^{n} a_{tj} x_j \right|
\le \sum_{j= 1,\, j \neq t}^{n} |a_{tj}|\ |x_j|.
\end{equation}
Applying the generalized H$\ddot{\mbox{o}}$lder inequality to (\ref{eqn9}), we have
\begin{eqnarray*}
 |\mu-a_{tt}||x_t| &\le& \left(\sum_{j= 1,\, j \neq t}^{n} |a_{tj}|^{p}\right)^{\frac{1}{p}}
 \left( \sum_{j= 1,\, j \neq t}^{n} |x_j|^{q}\right)^{\frac{1}{q}}.
\end{eqnarray*}
Since $|x_t| \ge |x_i|$  for all $1 \le i \le n$, we have
\begin{eqnarray*}
 |\mu-a_{tt}||x_t| &\le & n_t^{(p)}(A) \left( (n-1) |x_t|^{q} \right)^{\frac{1}{q}},
\end{eqnarray*}
 that is,
\begin{equation}\label{eq10}
 |\mu-a_{tt}| \le n_t^{(p)}(A) \left(n-1\right)^{\frac{1}{q}}.
\end{equation}
Similarly, using $|x_t| \ge |x_i|\,\,\forall\,\,i\,\, (1 \le i \le n)$ in (\ref{eqn9}), we get
\begin{equation}\label{eq11}
|\mu-a_{tt}|\le \sum_{j= 1,\, j \neq t}^{n} |a_{tj}|= r_t(A).
\end{equation}
Combining (\ref{eq10}) and (\ref{eq11}) for $\gamma \in [0,\,1],$ we have
\begin{equation}\label{eq12}
 |\mu- a_{tt}|^{1-\gamma}
\le (n_t^{(p)}(A))^{1-\gamma} (n-1)^{\frac{1-\gamma}{q}}\,\,\mbox{and}\,\,
|\mu-a_{tt}|^{\gamma}\leq r_t(A)^{\gamma},
\end{equation}
that is,
$$|\mu-a_{tt}| \le  (n-1)^{\frac{1-\gamma}{q}} (n_t^{(p)}(A))^{1-\gamma} r_t(A)^{\gamma}.\,\,\,\blacksquare$$

 Let us relate Theorem \ref{g1} to some existing results:
\begin{itemize}
\item  Setting $p= q= 2$ and $\gamma= 1$ implies that the left eigenvalues of $A:=(a_{ij}) \in M_n({\H})$
are  contained in the
 union of $n$ Greschgorin balls $B_i(A):= \left\{z\in \H : |z-a_{ii}| \le  r_i(A) \right\}, 1 \leq i \leq n,$ that is,
  \[ \Lam_l(A) \subseteq B(A):= \cup_{i=1}^{n} B_i(A).\]
This  result can be found in \cite[Theorem~6]{fz07}.
\item  Setting $p= q= 2$ and  $\gamma= 0$ implies that the
left eigenvalues of $A:=(a_{ij}) \in M_n({\H})$ are  contained in the
  union of $n$ balls $B_i(A):= \left\{z\in \H : |z-a_{ii}| \le
  (n-1)^{\frac{1}{2}} n_i^{(2)}(A) \right\}, 1 \leq i \leq n,$ that is,
 \[ \Lam_l(A) \subseteq B(A):=\cup_{i=1}^{n} B_i(A).\]
This result can be found in \cite[Theorem~1]{jw08}.
\end{itemize}

 We now present a generalization of \cite[Theorem~7]{fz07} and
 \cite[Theorem~3.1]{lyj12} by applying the generalized  H$\ddot{\mbox{o}}$lder
 inequality over the skew field of quaternions. For a general matrix $A:= (a_{ij}) \in M_n({\H})$ ,
 all the right eigenvalues may not lie in the union of $n$ generalized balls $B_i(A), 1 \leq i \leq n$.
 On the other hand,  we show that
 every connected region of the generalized balls $B_i(A), 1 \leq i \leq n$ contains
 some right eigenvalues of $A$.
 \begin{theorem}\label{gg2}
  Let $A:= (a_{ij}) \in M_n({\H})$ and let $\gamma \in [0, 1]$. For every
  right eigenvalue $\mu$ of $A$ there exists a nonzero quaternion $\beta$ such that
  $\beta^{-1} \mu \beta$ $($which is also a right eigenvalue$)$ is contained in the
  union of $n$ generalized balls
 \[B_i(A):= \left\{z\in \H : |z-a_{ii}| \le (n-1)^{\frac{1-\gamma}{q}}r_i(A)^{\gamma}
 (n_i^{(p)}(A))^{1-\gamma}\right\},\quad  1 \leq i \leq n,\]
 that is,
$$\left\{ z^{-1} \mu z : 0 \ne z \in \H\right\} \cap \cup_{i= 1}^{n} B_i(A)\ne\emptyset,$$
 where $p, q \in (1, \infty)$ with $\frac{1}{p}+\frac{1}{q}=1.$
\end{theorem}
\proof Let $\mu$ be a right eigenvalue of $A.$ Then there exists some nonzero vector
$x \in \H^n$ such that $A x = x \mu.$ Let $x:= [x_1, \ldots, x_n]^T\in \H^n$ and choose $x_t$ from $x$ as given in Theorem \ref{g1}.
Consider $\rho \in \H$ such that $x_t \mu= \rho x_t$. Then we have
\begin{equation}\label{9}
 |\rho-a_{tt}||x_t|= \left|\sum_{j=1,\, j \neq t}^{n} a_{tj} x_j \right|
\le \sum_{j=1,\, j \neq t}^{n} |a_{tj}|\ |x_j|.
\end{equation}
 Using the method from the proof of Theorem \ref{g1}, we have
\[
 |\rho-a_{tt}|
\le  (n-1)^{\frac{1-\gamma}{q}} (n_t^{(p)}(A))^{1-\gamma} r_t(A)^{\gamma}.\,\,\,\blacksquare
\]

Let us relate Theorem \ref{gg2} to some existing results:
\begin{itemize}
\item Substituting $p= q= 2$ and $\gamma= 1$, we obtain
$$ \{ z^{-1} \mu z : 0 \ne z \in \H \} \cap \cup_{i= 1}^{n} \{z\in \H : |z-a_{ii}| \le r_i(A) \} \ne \emptyset. $$
This result can be found in \cite[Theorem~7]{fz07}.
\item Substituting $p= q= 2$ and $\gamma= 0$, we get
\[ \{ z^{-1} \mu z : 0 \ne z \in \H \} \cap \cup_{i= 1}^{n}
 \left\{z\in \H : |z- a_{ii}| \le \sqrt{n-1} \ n_i^{(2)}(A) \right\} \ne \emptyset.\]
 This result can be found in \cite[Theorem~3.1]{lyj12}.
\end{itemize}

We next present a sufficient condition for the stability of a matrix $A\in M_{n}(\H).$
\begin{proposition}
 Let $A:= (a_{ij}) \in M_n({\H})$ and let $\gamma \in [0, 1]$. Assume that
 \begin{eqnarray}\label{SE1}
 \Re(a_{ii})+ (n-1)^{\frac{1- \gamma}{q}}{r_i}(A)^{\gamma} (n_i^{(p)}(A))^{1-\gamma} <0, \quad 1 \leq i \leq n,
  \end{eqnarray}
  where $\frac{1}{p}+ \frac{1}{q} = 1$ with $ p, q \in (1, \infty).$
 Then the matrix $A$ is stable.
\end{proposition}
\proof Let $\lam \in \Lam_r(A).$  From Theorem \ref{gg2} there
exists $0 \neq \rho \in \H$ such that $\rho^{-1} \lam \rho \in \cup _{i=1}^{n}B_i(A).$
Without loss of generality, we assume $\rho^{-1} \lam \rho \in B_l(A),$ that is,
$$
 |\rho^{-1} \lam \rho-a_{ll}| \leq (n-1)^{\frac{1- \gamma}{q}}{r_l}(A)^{\gamma} (n_l^{(p)}(A))^{1-\gamma}.
 $$
 Consider $\lam:=\lam_1+\lam_2 {\bf{i}}+\lam_3{\bf{j}}+\lam_4{\bf{k}}$ and $a_{ll}=a_l+b_l {\bf{i}}+c_l {\bf{j}}+d_l {\bf{k}}.$ Then
from (\ref{SE1}), we obtain
 \begin{eqnarray}\label{SE4}
 | (\lam_1 -a_l) +(\rho^{-1} \lam_2 {\bf{i}} \rho-b_l{\bf{i}} ) +(\rho^{-1} \lam_3 {\bf{j}}\rho-c_l{\bf{j}} )
 +(\rho^{-1} \lam_4 {\bf{k}}\rho-d_l{\bf{k}} )| < -\Re(a_{ll})=- a_l .
 \end{eqnarray}
 The equality (\ref{SE4}) is possible when $\lam_1 <0,$ that is, $\Re(\lam) <0,$ hence $\lam \in \H^{-}.$
This shows that the matrix $A$ is stable.$\,\,\,\blacksquare$

When all the diagonal entries of a matrix $A \in M_n(\H)$ are real, we have the following theorem.
  \begin{theorem}\label{t22}
   Let $A:= (a_{ij}) \in M_n({\H})$ with $a_{ii}\in \R$ and let $\gamma \in [0, 1]$. Then all the right eigenvalues of $A$
    are contained in the union of $n$ generalized balls $$B_i(A):= \left\{z\in \H : |z- a_{ii}| \le (n-1)^{\frac{1- \gamma}{q}}
    {r_i}(A)^{\gamma} (n_i^{(p)}(A))^{1-\gamma}\right\},\quad 1 \leq i \leq n,$$
    that is,
   $$ \Lambda_r(A)\subseteq B(A):= \cup_{i= 1}^{n}B_i(A),$$
   where $p, q \in (1, \infty)$ with  $\frac{1}{p}+\frac{1}{q}=1.$
  \end{theorem}
\proof
Let $\lam$ be a right eigenvalue of $A.$ Then there exists some nonzero vector $x \in \H^n$ such that $A x =x \lam$.
Let $x:=[x_1, \ldots, x_n]^T \in \H^n$ and let $x_t$ be an element of $x$ such that
$|x_t|\ge |x_i|, 1\le i \le n$. Then $|x_t|> 0$. Thus from $Ax= x \lam,$ we have
\begin{eqnarray*}
 a_{tt} x_t + \dm{\sum_{j=1, \,  j \neq t}^{n} a_{tj} x_j} &=&  x_t \lam,
\end{eqnarray*}
since $a_{tt} \in \R,$ so $a_{tt} x_t=x_t a_{tt}.$ Then from  the proof method of Theorem \ref{g1}, we have
 \[
 |\lam-a_{tt}|
\le  (n-1)^{\frac{1-\gamma}{q}} (n_t^{(p)}(A))^{1-\gamma} r_t(A)^{\gamma}.\,\,\,\blacksquare
\]

The above result has great significance as Hermitian and $\eta$-Hermitian matrices have all real diagonal entries. In general, $\eta$-Hermitian matrices arise widely in
applications \cite{cdf11, rf12, cd11}. To that end, we state the following proposition when all diagonal
entries of $A\in M_n (\H)$ are real.
In particular, this result gives a sufficient condition for the stability of a matrix $A \in M_n(\H)$.
\begin{proposition}\label{st1}
 Let $A:= (a_{ij}) \in M_n ({\H})$ with $a_{ii}\in \R$ and let $\gamma \in [0, 1]$. Assume that
 $$a_{ii}+ (n-1)^{\frac{1- \gamma}{q}}{r_i}(A)^{\gamma} (n_i^{(p)}(A))^{1-\gamma} < 0, \quad 1 \leq i \leq n,$$
 $ \mbox{ where}\,\, p, q \in (1,\,\infty)\,\,
 \mbox{with}\,\, \frac{1}{p}+ \frac{1}{q} = 1.$ Then the matrix $A$ is stable.
\end{proposition}


From Theorem \ref{t22}, all the complex right eigenvalues of a matrix
$A=(a_{ij}) \in M_n(\H)$ with all real diagonal entries lie in the union of $n$-discs
$D_i(A):= \{z \in \C: |z-a_{ii}| \leq (n-1)^{\frac{1- \gamma}{q}}
    {r_i}(A)^{\gamma} (n_i^{(p)}(A))^{1-\gamma} \}, 1 \leq i \leq n,$ that is,
    \begin{eqnarray}\label{EEE1}
   \Lambda_c(A)\subseteq D(A):= \cup_{i= 1}^{n}D_i(A).
   \end{eqnarray}
However, if diagonal entries are from $\C \setminus \R,$ then it is not necessary that  all the complex right eigenvalues of $A$ are contained
in the union of $n$-discs $D_i(A), 1 \leq i \leq n$ as the following examples suggest.
\begin{exam}\label{ree2}{\rm Let
$
 A:= \bmatrix{1-2{\bf{i}} & {\bf{j}} & {\bf{k}} \\ 0  & -2{\bf{i}} & -{\bf{i}} \\ 0 & {\bf{k}} & 3+{\bf{i}}}.
$
The set of complex right eigenvalues of $A$ is $$\Lambda_c(A):=\{\lam_1, \lam_2, \lam_3,\lam_4, \lam_5, \lam_6 \},$$
where $\lam_1= -0.0164+2.0083 {\bf{i}},\, \lam_2= -0.0164-2.0083 {\bf{i}},\, \lam_3=1+2{\bf{i}},\, \lam_4= 1-2{\bf{i}},\,
\lam_5= 3.0164+1.0324 {\bf{i}}$, and $ \lam_6=3.0164+1.0324 {\bf{i}}.$

For $\gamma= 1$ in (\ref{EEE1}), the discs $D_1(A), D_2(A),$ and $D_3(A)$ are as follows:
$$
 D_1(A):= \{z \in \C: |z-1+ 2{\bf{i}}| \leq 2 \}, \,\,\, D_2(A):= \{z \in \C: |z+2{\bf{i}}| \leq 1 \},\,\,\mbox{and}$$
 $$D_3(A):=\{z \in \C: |z-3-{\bf{i}}| \leq 1 \}.$$
 From Figure \ref{fig2}, it is clear that
$\lam_1, \lam_3,$ and $\lam_6$ lie outside the discs $D_1(A), D_2(A),$ and
$D_3(A).$
  \begin{figure}[h!]
   \begin{center}
   \includegraphics[height=7cm,width=7cm]{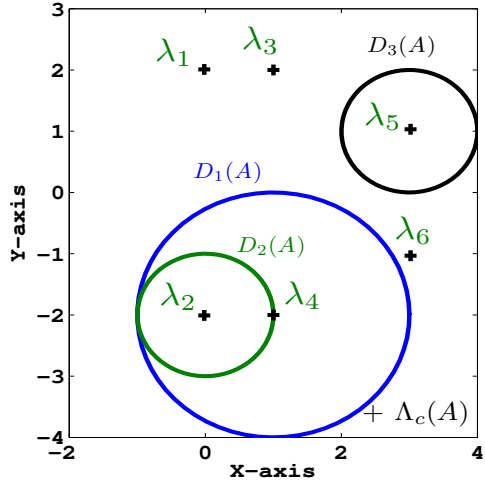}
    \caption{Location of the complex right eigenvalues of $A$ from Example \ref{ree2}.}\label{fig2}
   \end{center}
   \end{figure}
 }
\end{exam}
\begin{exam}\label{ree3}{\rm Let
$
 A=\bmatrix{-4 & 1+{\bf{j}}+\sqrt{2} {\bf{k}} & {\bf{j}}\\ {\bf{i+j}} & -10 & 2{\bf{j}}-{\bf{k}} \\ {\bf{i}}-2{\bf{j}}+2{\bf{k}} &
 \sqrt{3}+2{\bf{j}}-3{\bf{k}} & -8}.
$
In this example, there are six complex right eigenvalues $\lam_j\,\,( 1 \leq j \leq 6)$ which are 
shown in Figure \ref{fig3}.  Substituting $\gamma=1$ in (\ref{EEE1}), then all the complex right
eigenvalues of the matrix $A$ are contained in the union of three discs $D_1(A), D_2(A),$
and $D_3(A),$ where
$$
 D_1(A):=\{z \in \C: |z+4| \leq 3 \}, \,\,\, D_2(A):=\{z \in \C: |z+10| \leq \sqrt 2 +\sqrt 5\},\,\,\mbox{and}$$
 $$D_3(A):=\{z \in \C: |z+8| \leq 7 \}.$$
 From Figure \ref{fig3},  the standard right eigenvalues of $A$ are $\lam_1$, $\lam_3$, and $\lam_5$. Then
 \[
  \Lam_r(A)=[\lam_1]\cup[\lam_3]\cup[\lam_5].
 \]
 Also, from Figure \ref{fig3}, we observe that $\Re(\lam_i) \in \H^{-}\, (i=1,3,5).$ Hence
 \[
  \Re(\lam_1)=\Re(\rho^{-1} \lam_1 \rho),
  \,\Re(\lam_2)=\Re(\tau^{-1} \lam_2 \tau), \mbox{and}\,\Re(\lam_3)=\Re(\nu^{-1} \lam_3 \nu)\,\,\,\, \forall \,
  \rho, \tau, \nu \in \H
 \]
 Thus the matrix $A$ is stable.

  \begin{figure}[h!]
   \begin{center}
   \includegraphics[height=7cm,width=7cm]{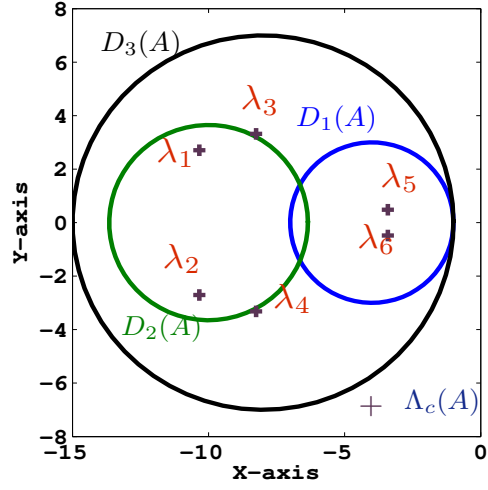}
   \caption{Location of the complex right eigenvalues of $A$ from Example \ref{ree3}.}\label{fig3}
   \end{center}
   \end{figure}

 }
\end{exam}

 In general, similar quaternionic matrices may not have the same left eigenvalues, see, \cite[Example 3.3]{fz07}. However, the following result is true.
 \begin{proposition}\label{sml}
  Let $A \in M_n (\H)$ and let $W$ be any invertible real matrix. Then $A$ and $WAW^{-1}$ have the same left eigenvalues.
 \end{proposition}
\proof Let $\lam$ be a left eigenvalue of $A.$ Then there exists some nonzero vector $x\in \H^n$
such that $Ax= \lam x.$ Let $W$ be an invertible real matrix. Then
\[
 WAx= W\lam x= \lam Wx.
\]
Now,  $WAW^{-1} Wx=\lam Wx.$ Setting $Wx= y$  implies $W A W^{-1} y=\lam y.\,\,\blacksquare$

Let $A:=(a_{ij}) \in M_n(\H).$ Suppose $W=\diag(w_1, w_2, \ldots, w_n)$ with $w_i \in \R^+, 1 \leq i \leq n.$ Then
\[
 W^{-1} A W=\left(\frac{a_{ij} w_j}{w_i}\right)\,\,\mbox{and}\,\, \Lam_l(A)= \Lam_l(W^{-1} A W).
\]
Define
$$ r_i^W(A):= \sum_{j=1,\, j\neq i}^{n} \frac{|a_{ij}| w_j}{w_i}\,\,\,\mbox{and}\,\,\,
\dm{ c_i^W(A):= \sum_{j=1,\, j\neq i}^{n} \frac{|a_{ji}| w_i}{w_j}},\quad  1 \leq i \leq n.$$

Applying Theorem \ref{os1} to $W^{-1} A W$, we get the following theorem which may be sharper than Theorem \ref{os1} depending
upon the choice of $W$.
\begin{theorem}\label{osts}
 Let $A:= (a_{ij}) \in M_n(\H)$. Then all the left eigenvalues of $A$ are contained in the union of $n$
 balls 
 $$T_i^W(A):= \{z \in \H : |z-a_{ii}| \leq (r_i^W(A))^{\gamma}\, (c_i^W(A))^{1-\gamma}\},\quad  1 \leq i \leq n,$$
 that is,
 \[
 \Lam_l(A)= \Lam_l(W^{-1} A W) \subseteq T^W(A):=\cup_{i=1}^{n} T_i^W(A).
\]
%
%
\end{theorem}
Since the above theorem holds for every $W=\diag(w_1, w_2, \ldots, w_n)$, where $w_i \in \R^+,$ we have
\[
 \Lam_l(A)=\Lam_l(W^{-1} A W) \subseteq \underset{\substack{W \in M_n(S)}}{\cap} T^W(A)=: T^S(A),
\]
where $M_n(S)$ is a set of real diagonal matrices with non-negative entries. $T^S(A)$ is called the
minimal Ostrowski type set for the matrix $A$.

Substituting $\gamma= 1$ in Theorem \ref{osts}, we obtain
\begin{equation}\label{neqn1}
 \Lam_l(A)=\Lam_l(W^{-1} A W) \subseteq \eta^W (A):=\cup_{i=1}^{n} \eta_i^W(A),
\end{equation}
where $ \eta_i^W(A):=\left\{z \in \H : |z-a_{ii}| \leq r_i^W(A) \right\}.$
 Therefore,
\[
 \Lam_l(A)=\Lam_l(W^{-1} A W) \subseteq \underset{\substack{W \in M_n(S)}}{\cap} \eta^W(A)=: \eta^S(A),
\]
where $\eta^S(A)$ is called the first minimal Gerschgorin type set for the matrix $A$.

For $\gamma= 0$ in Theorem \ref{osts}, we have
\begin{eqnarray}\label{neqn2}
 \Lam_l(A)=\Lam_l(W^{-1} A W) \subseteq \Omega^W(A):=\cup_{i=1}^{n} \Omega_i^W(A),
 \end{eqnarray}
 where
 $
\Omega_i^W(A):=\left\{z \in \H : |z-a_{ii}| \leq c_i^W(A) \right\}.
$
Then
\[
 \Lam_l(A)=\Lam_l(W^{-1} A W) \subseteq \underset{\substack{W \in M_n(S)}}{\cap} \Omega^W(A)=: \Omega^S(A),
\]
where $\Omega^S(A)$ is called the second minimal Gerschgorin type set for the matrix $A$.

Equivalently, applying Theorem \ref{rc} to $ W^{-1} A W,$ we get the following theorem:

\begin{theorem}\label{octs}
 Let $A:=(a_{ij}) \in M_n (\H)$ and let $\gamma \in [0, 1]$. Then all the left eigenvalues of $A$ are contained in the union of $\frac{n(n-1)}{2}$
 ovals of Cassini
 $$K_{ij}^W(A):=\{z \in \H : |z-a_{ii}|\ |z-a_{jj}|\leq \left(r_i^W(A) \right)^{\gamma} (r_j^W(A))^{\gamma}
(c_i^W(A))^{1-\gamma} (c_j^W(A))^{1-\gamma} \},\quad 1 \leq i, j \leq n,\quad i \neq j,$$
that is,
\[
  \Lam_l(A)=\Lam_l(W^{-1} A W) \subseteq K^W(A):=\cup_{\substack{i, j=1\\ i \neq j}}^{n} K_{ij}^W (A).
  \]
\end{theorem}
Since Theorem \ref{octs} holds for every $W= \diag(w_1, w_2, \ldots, w_n)$ with $w_i \in \R^+.$ Then
\[
 \Lam_l(A)= \Lam_l(W^{-1} A W) \subseteq \underset{\substack{W \in M_n(S)}}{\cap} K^W(A)=: K^S(A).
\]
$K^S(A)$ is called the minimal Brauer type set for the matrix $A$.
\begin{exam}\label{se}{\rm
 Let
 $
  A=\bmatrix{ {\bf{j}} &   {\bf{k}} & {\bf{2j+\sqrt{5}k}} \\0 & {\bf{i+k}} &      {\bf{\sqrt{2}i+j-k}}     \\ 0 & 0 & {\bf{2-i}}   }.
  $
  Let $\gamma= 1$ in Theorem \ref{os1}. Then, we have
the three Gerschgorin type balls $G_1(A):=\{ z \in \H : |z-{\bf{j}}| \leq 4 \}, G_2(A):=\{ z \in \H : |z-{\bf{i-k}}| \leq 2\},$ and
 $ G_3(A):=\{ z \in \H : |z-{\bf{2+i}}| \leq 0\}.$
 If $W=\diag(w_1, w_2, w_3)$ with $w_1=8,\ w_2=4, w_3=1.$ Then by (\ref{neqn1})
 $$\eta_1^W(A):= \{ z \in \H : |z-{\bf{j}}| \leq 7/8\},\,\,\eta_2^W(A):= \{ z \in \H : |z-{\bf{i-k}}| \leq 1/2\},\,\,\mbox{and}$$
 $$\eta_3^W(A):= \{ z \in \H : |z-{\bf{2+i}}| \leq 0\}.$$
 Hence it is clear that $\eta_1^W(A) \subset G_1(A)$ and $\eta_2^W(A) \subset G_2(A).$

 For $\gamma=1$, Theorem \ref{rc} gives
 the following ovals of Cassini:
 \[
  K_{12}(A):= \{ z \in \H : |z-{\bf{j}}|\, |z-{\bf{i-k}}| \leq 8 \},\
 K_{23}A):= \{ z \in \H :  |z-{\bf{i-k}}|\, |z-{\bf{2+i}}| \leq 0 \},\,\mbox{and}
 \]
 \[
  K_{31}(A):=\{ z \in \H : |z-{\bf{2+i}}| \,  |z-{\bf{j}}| \leq 0 \}.
  \]
Consider $W=\diag(w_1, w_2, w_3)$ with $w_1=w_2=6,$ and $\,w_3= 1.$ Then by Theorem \ref{octs} with $\gamma= 1,$ we obtain
 \[
K_{12}^W(A):= \{ z \in \H : |z-{\bf{j}}|\, |z-{\bf{i-k}}| \leq 1/2 \},\,\,
  K_{23}^W(A):= \{ z \in \H :  |z-{\bf{i-k}}|\, |z-{\bf{2+i}}| \leq 0 \},\,\,\mbox{and}
  \]
  \[
 K_{31}^W(A):= \{ z \in \H : |z-{\bf{2+ i}}| \,  |z-{\bf{j}}| \leq 0 \}.
 \]
   Hence $K_{12}^W(A) \subset K_{12}(A).\,\,\,\blacksquare$
   }
\end{exam}
\section{Bounds for the zeros of quaternionic polynomials}\label{s4}
 In this section, we derive bounds for the zeros of quaternionic polynomials by applying the 
localization theorems for the left eigenvalues of a quaternionic matrix. Due to noncommutivity of quaternions, we first define
some basic facts on multiplication
of quaternions.
For $p, q \in \H$, define $p \times q:= pq.$ For $0 \neq p\in \H$ and $q \in \H$, define
\[
\frac{1}{p} \times q:= p^{-1}\times q:= p^{-1}q,\,q \times \frac{1}{p}:= q \times p^{-1}:= qp^{-1}.
\]
Recall the quaternionic polynomials $p_l(z)$ and $p_r(z)$ from (\ref{l}) and (\ref{lk}).
Then the corresponding companion matrices of the simple monic polynomials $p_l(z)$ and $p_r(z)$ are given by
\[
C_{p_l}:= \bmatrix{0 & \vrule & 1 & & 0\\ \vdots & \vrule & & \ddots &\\0 & \vrule & 0 & &1\\
 \cline{1-5}
-q_0 & \vrule & -q_1 &\ldots & -q_{m-1}}:= \kbordermatrix{ &1& & m-1\\
 m-1& 0 & \vrule & I \\
  \cline{2-4}
 1& C_{p_l}(m, 1) & \vrule   &  C_{p_l}(m, 2:m)
 }\, \mbox{and}\, C_{p_r}:=C_{p_l}^T,
 \]
respectively.
%
 Let $q_0 \neq 0$, and define simple monic reversal polynomials of $p_l(z)$ and $p_r(z)$ as follows:
\[q_l(z):= \frac{1}{q_0} \times p_l\left(\frac{1}{z}\right) \times z^m= z^m+ q_0^{-1} q_1 z^{m-1}
+\dots+ q_0^{-1}q_{m-1}z+ q_0^{-1},\]

\[ q_r(z):= z^m  \times p_r\left(\frac{1}{z}\right) \times \frac{1}{q_0} = z^m+ z^{m-1}q_1 q_0^{-1}+
 \dots+zq_{m-1}q_0^{-1} +q_0^{-1},
 \]
 respectively. 
 The corresponding companion matrices of the simple monic reversal polynomials $q_l(z)$ and $q_r(z)$ are denoted by $C_{q_l}$ and
 $C_{q_r}$, respectively.
 We observe that the zeros of $q_l(z)$ and $q_r(z)$ are the reciprocal of zeros of $p_l(z)$ and $p_r(z),$
 respectively.

 Now, we need the following result:
 \begin{proposition}\cite[Proposition 1]{rp01}.\label{lr2}
Let  $\lam \in \H.$ Then $\lam$ is a zero of the simple monic polynomial $p_l(z)$  if and only if $\lam$ is a left eigenvalue of its
corresponding companion matrix $C_{p_l}$.
 \end{proposition}

In general, a right eigenvalue of $C_{p_l}$ is not necessarily a zero of the simple monic polynomial $p_l(z)$.
For example, let a simple monic polynomial $p_l(z)= z^2+ {\bf{j}}z+ 2.$ Then its companion matrix is given by
 $$C_{p_l}=\bmatrix{0 & 1\\-2 &-{\bf{j}}}.$$ Here ${\bf{i}}$ is a right eigenvalue of $C_{p_l}.$ However, ${\bf{i}}$ is not a
 zero of $p_l(z).$

Analogous to Proposition \ref{lr2}, the following result is presented for $p_r(z)$.
\begin{proposition}\label{lr22}
Let  $\lam \in \H.$ Then $\lam$ is a zero of the simple monic polynomial $p_r(z)$  if and only if $\lam$ is a left eigenvalue of its
corresponding companion matrix $C_{p_r}.$
 \end{proposition}

We now present bounds for the zeros of $p_l (z)$ as follows.
%
%
\begin{theorem}\label{u1}
 Let $p_l(z)$ be a simple monic polynomial over $\H$ of degree $m.$ Then every zero $\tilde{z}$ of $p_l(z)$
 satisfies the following inequality:
$$\left(\dm{\max_{1 \leq i \leq m}}\left( r_i'(C_{q_l})^{\gamma}\, c_i'(C_{q_l})^{1-\gamma}\right)\right)^{-1} \leq |\tilde{z}|
\leq \dm{\max_{1 \leq i \leq m}}\left( r_i'(C_{p_l})^{\gamma}\, c_i'(C_{p_l})^{1-\gamma}\right),
$$
for every $\gamma \in [0, 1].$
\end{theorem}
\proof
From Proposition \ref{lr2},  zeros of $p_l(z)$ and left eigenvalues of
$C_{p_l}$ are same.
Thus, if
$\tilde{z}$ is a zero of $p_l(z),$ then $\tilde{z}$ is a left eigenvalue of $C_{p_l}.$
By applying Theorem \ref{os1} (Ostrowski type theorem) to $C_{p_l}$, we obtain
\[
|\tilde{z}|\leq \dm{\max_{1 \leq i \leq m}}\left( r_i'(C_{p_l})^{\gamma}\, c_i'(C_{p_l})^{1-\gamma}\right).
\]
We use the respective upper bounds for the zeros of the simple monic reversal polynomial $q_l(z)$ for the desired lower bounds for the
 zeros of $p_l(z)$. $\,\,\, \blacksquare$
\begin{corollary}\label{co}
Let $p_l(z)$ be a simple monic polynomial over $\H$ of degree $m.$ Then every zero $\tilde{z}$ of $p_l(z)$
 satisfies the following inequalities:
 \begin{enumerate}
  \item
$\dm{\frac{|q_0|}{\dm{\max_{1 \leq i \leq (m-1)}}\left\{1, |q_0|+ |q_i|\right\}}} \leq |\tilde{z}|
  \leq \dm{\max_{1 \leq i \leq (m-1)}} \left\{|q_0|, 1+|q_i|\right\}.$
\item
$  \dm{\frac{|q_0|}{\max\left\{|q_0|,1+\sum_{i=1}^{m-1} |q_i|\right\}}}\leq |\tilde{z}| \leq
\max\left\{1, \sum_{i=0}^{m-1} |q_i|\right\}. $
 \end{enumerate}
\end{corollary}
\proof Substituting $\gamma=0, 1$ in Theorem \ref{u1}, we obtain the desired results.$\,\,\,\blacksquare$

Next, we derive the following lemma which gives a better bound than Opfer's bound \cite[Theorem 4.2]{go09} for $|q_0|\geq 1$.
\begin{lemma}\label{l1}
 Assume that $|q_0| \geq 1.$ Then $
 \alpha \leq \mathcal{T},$
 where $\alpha:=\dm{\max_{1 \leq i \leq m-1}} \left\{|q_0|, 1+|q_i|\right\}\, \mbox{and}\,\,
  \mathcal{T}:=\max\left\{1, \sum_{i=0}^{m-1} |q_i|\right\}.$
\end{lemma}
\proof  {\bf{Case 1:}} If $|q_0|=1,$ then
\[
 \alpha=\dm{\max_{1 \leq i \leq m-1}} \left\{|q_0|, 1+|q_i|\right\}
 =\dm{\max_{1 \leq i \leq m-1}} \left\{ 1+|q_i|\right\}. \,\,\mbox{Also}
 \]

 $\mathcal{T}:=\max\left\{1, \sum_{i=0}^{m-1} |q_i|\right\}=\max\left\{1, |q_0|+\sum_{i=1}^{m-1} |q_i|\right\}=
1+\sum_{i=1}^{m-1} |q_i|.
$

{\bf{Case 2:}} If $|q_0| >1,$ then
\[
 \alpha=\dm{\max_{1 \leq i \leq (m-1)}} \left\{|q_0|, 1+|q_i|\right\}= |q_0| \,\,\mbox{or}\,\, \max_{1 \leq i \leq (m-1)}\left\{
 1+|q_i|\right\}\,  \mbox{and}
 \]
 $\mathcal{T}:=\max \{1, \sum_{i=0}^{m-1} |q_i| \}=\max\left\{1, |q_0|+ \sum_{i=1}^{m-1} |q_i|\right\}=
 |q_0|+\sum_{i=1}^{m-1} |q_i|.$
 Thus $\alpha \leq \mathcal{T}.$ This completes the proof. $\,\,\,\blacksquare$

On the other hand, if $|q_0|<1,$ then $\alpha \leq \mathcal{T}$ or $\alpha > \mathcal{T}$.  For example, for a simple
monic polynomial
$
 p'_l(z):= z^3+({\bf{i}}+2{\bf{j}}+2{\bf{k}})z^2-2{\bf{k}}z+0.5 {\bf{k}},
$
we have $\alpha=4$ and $\mathcal{T}=5.5.$ Hence $\alpha < \mathcal{T}$. Further, if we consider
$
 p''_l(z)= z^3+0.5{\bf{j}}z^2+(0.2{\bf{i}}+0.3{\bf{j}})z+0.5 {\bf{i}},
$
then $\alpha= 1.5$ and $\mathcal{T}= 1.36.$ Hence $\alpha > \mathcal{T}.$

 Next, by applying Theorem \ref{os1} to $WC_{p_l}W^{-1}$ and $WC_{q_l}W^{-1}$ ($W$ is an invertible real diagonal matrix), we
obtain different and potentially sharper bounds.
\begin{theorem}\label{sm}
     Let $w_i \in \R^+$, $1\le i \le m.$ Then every zero $\tilde{z}$ of the simple monic polynomial $p_l(z)$ satisfies
    the following inequality:
      $$ \left[\dm{\max_{1 \leq i \leq m}}\left\{ r_i'(W C_{q_l}W^{-1})^{\gamma}\,
c_i'(WC_{q_l}W^{-1})^{1-\gamma}\right\}\right]^{-1} \leq |\tilde{z}|
\leq \dm{\max_{1 \leq i \leq m}}\left\{ r_i'(WC_{p_l}W^{-1})^{\gamma}\, c_i'(WC_{p_l}W^{-1})^{1-\gamma}\right\},$$
 where $W:=\diag(w_1, w_2, \ldots, w_m)$ and $\gamma \in [0, 1].$
\end{theorem}
\proof  The companion matrix of $p_l(z)$ is given by
\[C_{p_l}=\kbordermatrix{ & 1 & & m-1\\
 m-1& 0 & \vrule & I \\
  \cline{2-4}
 1&-q_0  & \vrule & [-q_1 \ldots -q_{m-1}]
 }.\]
 Then
 \[
W C_{p_l} W^{-1}=\kbordermatrix{ &1& & m-1\\
 m-1& 0 &\vrule & \diag\left(\frac{w_1}{w_2}, \ldots, \frac{w_{m-1}}{w_m} \right) \\
  \cline{2-4}
 1&- \frac{w_m }{w_1} q_0  & \vrule & -\frac{w_m }{w_2} q_1 \ldots  -q_{m-1}
 }.
\]
By Proposition \ref{sml}, $C_{p_l}$ and $W C_{p_l} W^{-1}$
 have the same left eigenvalues. Rest of the proof follows from the proof method
 of Theorem \ref{u1}.$\,\,\,\blacksquare$
\begin{corollary}\label{cs}
      Let $p_l(z)$ be a simple monic polynomial over $\H$ of degree $m.$ Then every zero $\tilde{z}$ of $p_l(z)$
      satisfies the following inequalities:
\begin{enumerate}
\item $
  \left[\dm{\max_{0\le j \le m-1}\left\{\frac{(|q_0|w_j+ w_m|q_{m-j}|)}{|q_0|d_{j+1}}\right\}} \right]^{-1}\le
  \dm{|\tilde z| \le \dm{\max_{0\le j \le m-1}} \left\{\frac{w_j+ w_m|q_j|}{w_{j+1}}\right\}},\,\mbox{ where}\,\,\, w_0= 0.
 $

\item $  \left[\dm{\max_{1 \leq j\leq m-1}\left\{\frac{w_j}{w_{j+1}},\sum_{i=0}^{m-1}\frac{w_m|q_i|}{|q_0|w_{i+1}}\right\}}\right]^{-1}\leq
|\tilde{z}| \leq \dm{\max_{1 \leq j\leq m-1}\left\{\frac{w_j}{w_{j+1}},\sum_{i=0}^{m-1}\frac{w_m|q_i|}{w_{i+1}}\right\}}.$
 \end{enumerate}
\end{corollary}
\proof Substituting $\gamma= 0, 1$ in Theorem \ref{sm}, we get the desired results.$\,\,\,\blacksquare$

Let $w_j= w_m |q_j|, 1 \leq  j \le m-1,$ in the part (1) of Corollary \ref{cs}. Then we obtain
\[
 |\tilde z| \le \max_{1\le j \le m-1} \left\{\left|\frac{q_0}{q_1}\right|, 2 \left| \frac{q_j}{q_{j+1}} \right| \right\}.
\]
This is called the Kojima type bound for the zeros of the simple monic polynomial $p_l(z).$

 %
%
%
%
For computation of bounds of the zeros of $p_r(z)$, we define the following polynomial:
\[\tilde{p_l}(z):= \overline{p_r(\overline{z})}:= \sum_{j=0}^{m}  \overline{q_j} z^j.\]

Now, we discuss the following theorem which shows relation between the zeros of $p_r(z)$ and $\tilde{p_l}(z).$
\begin{theorem}\label{conj} Let $\lam \in \H.$ Then $\lam$ is a zero of the simple monic
polynomial $p_r(z)$ if and only if $\overline{\lam}$ is a zero of the simple monic polynomial
 $\tilde{p_l}(z).$
\end{theorem}
\proof The corresponding companion matrices of $p_r(z)$ and $\tilde{p_l}(z)$
are given by
$$
 C_{p_r}:=C_{p_l}^T\,\,\mbox{and}\,\, C_{\tilde{p_l}}:= C^H_{p_r},$$
respectively. By Lemma \ref{prop4}, if $\lam$ is a left eigenvalue of $C_{p_r}$, then $\overline{\lam}$ is a left eigenvalue of
	$C^H_{p_r}=C_{\tilde{p_l}}.$ By Propositions \ref{lr2} and \ref{lr22}, the left eigenvalues of $C_{p_r}$ and $C_{\tilde{p_l}}$ imply the zeros of
$p_{r}(z)$ and $\tilde p_{ l}(z),$ respectively. Hence if  $\lam$ is a zero of
$p_r(z),$ then $\overline{\lam}$ is also a zero of
${\tilde{p_l}} (z).\,\,\,\blacksquare$
\begin{remark}
Similar results can be obtained for the quaternionic polynomial $p_r(z)$ as well.
\end{remark}

\section{Bounds for the zeros of quaternionic polynomials by using the powers of companion matrices}\label{s6}
First, we present some preliminary
results for the powers of companion matrices $C_{p_l}$ and $C_{p_r}.$ In general, if $\lam$ is a left eigenvalue of
a quaternionic matrix $A,$ then $\lam^2$ is not necessarily a
left eigenvalue of $A^2$. For example, for a quaternionic matrix $A= \bmatrix{0 & {\bf{i}} \\ -{\bf{i}} & 0},$ we have
$
  \Lam_l(A):=\left\{ \mu: \mu= \alpha+\beta {\bf{j}}+ \gamma{\bf{k}}, \alpha^2+\beta^2+ \gamma^2= 1\right\}$ and
  $A^2= \bmatrix{1 & 0 \\ 0 & 1}.$ So $\Lam_l(A^2):= \{1\}.$  Here ${\bf{j}}$ is a left eigenvalue of $A$  but ${\bf{j}}^2$ is not
 a left eigenvalue of $A^2.$

Now we prove the following result for left eigenvalues of $C_{p_l}$ and $C^t_{p_l}$ ($t$ is a nonzero integer).
\begin{proposition}\label{p2}
If $\lam$ is a left eigenvalue of $C_{p_l}$ with respect to the eigenvector $x \in \H^{n}$, then
  $\lam^{t}$ is a left eigenvalue of $C_{p_l}^{t}$
  corresponding to the same eigenvector $x \in \H^{n}$.
\end{proposition}
\proof {\bf Case (a)}: Let $t$ be a positive integer and let $\lam$ be a left eigenvalue of
 $C_{p_l}$. Then, there
 exists $ 0 \neq x:=\left[1, \lam, \lam^2, \ldots, \lam^{m-1}\right]^T \in \H^n$ such that $C_{p_l} x=\lam x.$ Therefore,
 \begin{eqnarray*}
C^2_{p_l} x&=&C_{p_l}(C_{p_l} x)=C_{p_l} x \lam=x\lam^2\\
\vdots\\
C^t_{p_l} x &=& C^{t-1}_{p_l} (C_{p_l}x)=C^{t-1}_{p_l} x \lam=\dots=x\lam^t =  \lam^t x.
\end{eqnarray*}
 Thus, $\lam^t$ is a left eigenvalue of matrix $C^t_{p_l}$ corresponding to the same eigenvector $x \in \H^{n}.$
 
 {\bf Case (b)}:   Let $t$ be a negative integer. From {\bf Case (a)}, we have
 $C_{p_l} x=x \lam $. This implies $C_{p_l}^{-1} x=x \lam^{-1}$.
 Therefore,
 \begin{eqnarray*}
C^{-2}_{p_l} x&=&C_{p_l}^{-1}(C_{p_l}^{-1} x)=C_{p_l}^{-1} x \lam^{-1}=x\lam^{-2}\\
\vdots\\
C^{t}_{p_l} x &=& C^{(t+1)}_{p_l} (C_{p_l}^{-1}x)=C^{(t+1)}_{p_l} x \lam^{-1}=\dots=x\lam^{t} =  \lam^{t} x.
\end{eqnarray*}
 Thus, $\lam^{t}$ is a left eigenvalue of $C^{t}_{p_l}$ with respect to the same eigenvector $x \in \H^{n}.\,\,\, \blacksquare$

Next, we state the following result for left eigenvalues of $C_{p_r}$ and $C^t_{p_r}$ ($t$ is a nonzero integer).
 \begin{proposition}\label{p3}
 If $\lam$ is a left eigenvalue of $C_{p_r}$ with respect to the eigenvector $x \in \H^{n}$, then
  $\lam^{t}$ $(t\,\, \mbox{is a nonzero integer})$ is a left eigenvalue of $C_{p_r}^{t}$
  corresponding to the same eigenvector $x \in \H^{n}$.
\end{proposition}
\proof {\bf Case (a)}:  Let $t$ be a positive integer and let $\lam$ be a left eigenvalue of
$C_{p_r}.$ Now from Lemma \ref{prop4}, $\overline{\lam}$ is a
left eigenvalue
of $C^H_{p_r}.$ Then there exists $0 \neq x:= \left[1, \overline{\lam}, (\overline{\lam})^2, \ldots, (\overline{\lam})^{m-1}\right]
\in \H^n$
such that $C^H_{p_r} x= \overline{\lam} x= x \overline{\lam}$. This gives
 \begin{eqnarray*}
\left(C^H_{p_r}\right)^2 x &=& C^H_{p_r}(C_{p_r}^H x)= C^H_{p_r} x \overline{\lam}= x (\overline{\lam})^2\\
\vdots\\
\left(C^H_{p_r}\right)^t x &=& \left(C^H_{p_r}\right)^{t-1} (C^H_{p_r}x)=\left(C^H_{p_r}\right)^{t-1} x \overline{\lam}=\dots=
x (\overline{\lam})^t =(\overline{\lam})^t x.
\end{eqnarray*}
Thus, $(\overline{\lam})^t$ is a left eigenvalue of $\left(C^H_{p_r}\right)^t.$ Then by Lemma \ref{prop4}, $\lam^t$ is a left eigenvalue of
$C^t_{p_r}.$

{\bf Case (b)}:  Let $t$ be a negative integer. From {\bf Case (a)}, we have
$C^H_{p_r} x= \overline{\lam} x= x \overline{\lam}$. This implies $(C^H_{p_r})^{-1} x=x (\overline{\lam})^{-1}$. Thus
 \begin{eqnarray*}
(C^{H}_{p_r})^{-2} x&=&(C^H_{p_r})^{-1}\{(C^H_{p_r})^{-1} x\}=(C^H_{p_r})^{-1} x (\overline{\lam})^{-1}=x (\overline{\lam})^{-2}\\
\vdots\\
 (C^{H}_{p_r})^t x &=& (C^H_{p_r})^{(t+1)} \{(C^H_{p_r})^{-1}x\}=(C^H_{p_r})^{(t+1)} x (\overline{\lam})^{-1}=\dots=
 x (\overline{\lam})^{t} =(\overline{\lam})^{t} x. 
\end{eqnarray*}
Thus, $(\overline{\lam})^t$ is a left eigenvalue of $\left(C^H_{p_r}\right)^t.$ Then by Lemma \ref{prop4}, $\lam^t$ is a left eigenvalue of
$C^t_{p_r}.\,\,\,\blacksquare$

Further, we present a framework to find the powers of the companion matrix $C_{p_l}$ which can be derived in a simple
procedure as follows, keeping in view that quaternions do not commute.
\begin{theorem}\label{i} Consider
$C_{p_l}=\kbordermatrix{ &1& & m-1\\
 m-1& 0 &\vrule & I \\
  \cline{2-4}
 1&C_{p_l}(m, 1)  & \vrule & C_{p_l} (m, 2: m)
 }$. \\
 {\bf (a)} If $t < m$ is a positive integer, then 
 \begin{equation}\label{q11}
 C_{p_l}^t=
 \kbordermatrix{ &t& & m-t\\
 m-t&  0 &\vrule & I \\
  \cline{2-4}
 t& C  & \vrule & D
 },
 \end{equation}
 {\bf (b)} if $t \ge m,$ then
%
 \begin{eqnarray}\label{q2}
  C^t_{p_l}= \left[
 \begin{array}{cc}
 C^{t-(m-1)}_{p_l}(m,1: m)\\
 C^{t-(m-2)}_{p_l}(m,1: m)\\
 \vdots\\
 C^{t-1}_{p_l}(m,1: m)\\
 C_{p_l}^t(m, 1:m)
 \end{array}\right]_{m\times m},
 \end{eqnarray}
where
\begin{eqnarray*}
C_{p_l}^t(m, 1) &:=& C^{t-1}_{p_l}(m, m) C_{p_l}(m,1),\,\,\\  C_{p_l}^t(m, 2: m) &:=& C_{p_l}^{t-1}(m, 1: m-1)+C_{p_l}^{t-1}(m, m)C_{p_l}(m, 2: m),
\end{eqnarray*}
$$C:= \bmatrix{C_{p_l}(m,1:t)\\C^2_{p_l}(m,1:t)\\
\vdots\\
C^t_{p_l}(m,1:t)}_{t\times t},\,\, \mbox{and}\,\, \, D:= \bmatrix{C_{p_l}(m,t+1:m)\\C^2_{p_l}(m,t+1:m)\\
\vdots\\
C^t_{p_l}(m,t+1:m)}_{t\times(m-t)}.$$
Note that $C_{p_l}(k, 1: m)$ denotes the $k$-th row of the matrix $C_{p_l}.$
\end{theorem}
\proof Assuming $t=1$, (\ref{q11}) becomes
 $C_{p_l}=\kbordermatrix{ &1& & m-1\\
 m-1& 0 & \vrule & I \\
  \cline{2-4}
 1& C_{p_l}(m, 1) & \vrule   &  C_{p_l}(m, 2: m)
 },$ where $C_{p_l}(m, 1):=-q_0,  C_{p_l}(m, 2: m):=[-q_1 \ldots -q_{m-1}].$ Thus the theorem is true for $t=1.$
 Now, let us consider $C_{p_l}$ as
\[
 C_{p_l}=
 \kbordermatrix{ &m-k& & k\\
  k &  A' &\vrule & B' \\
  \cline{2-4}
 m-k & C'  & \vrule & D'
 },\,\, \mbox{where}
\]
$A' := C_{p_l}(1: k, 1: m-k), B':= C_{p_l}(k+1: m, m-k+1: m), C':= C_{p_l}(k+1: m, 1: m- k), D': =C_{p_l}(k+1: m, m-k+1: m).$
For $t= k= 3$, we get
\begin{eqnarray*}
C_{p_l}^3&=&\kbordermatrix{ &2 & & m-2\\
 m-2&  0 &\vrule & I \\
  \cline{2-4}
 2 & C  & \vrule & D
 } \kbordermatrix{ &m-2& & 2\\
  2 &  A' &\vrule & B' \\
  \cline{2-4}
 m-2 & C'  & \vrule & D'
 }= \kbordermatrix{ &m-2& & 2\\
 m-2& C' &\vrule & D' \\
  \cline{2-4}
 2 & CA'+ DC'  & \vrule & CB'+ DD'
 }.
 \end{eqnarray*}
 Note that in each step, size of the identity matrix $I$ reduces by order $1$ and the size of matrix $C$
 increases by order $1.$ Similarly, the matrix $D$ increases by $1$ row and decreases by $1$ column. Finally,
 after rearranging and separating $0$ and $I$ matrices we get
$$
\kbordermatrix{ &2+1& & m-2-1\\
 m-2-1& 0 &\vrule & I \\
  \cline{2-4}
 2+1 & C  & \vrule & D
 },$$ where $C$ and $D$ are of size $3\times 3$ and $3\times (m-3),$ respectively. Assuming that the theorem is true for $t=k$, we have
\begin{eqnarray*}
 C_{p_l}^{k+1}= C_{p_l}^k C_{p_l}&=& \kbordermatrix{ &m-k& & k\\
 m-k& C' &\vrule & D' \\
  \cline{2-4}
 k & CA'+DC'  & \vrule &CB'+DD'
 }
=
\kbordermatrix{ &k+1& & m-k-1\\
 m-k-1& 0 &\vrule & I \\
  \cline{2-4}
 k+1 & C  & \vrule & D
 },
\end{eqnarray*}
where the corresponding $C$ and $D$ matrices are given in the statement of the theorem.

\noin The proof for $t\ge m$ is similar$.\,\,\,\blacksquare$

In the case of quaternionic matrix,  $C_{p_l}= C_{p_r}^T$ but
$C_{p_r}^t \not= (C_{p_l}^t)^T $  for $t \ge 2.$  This is illustrated by the following example.

 \begin{exam}\label{ex5.7}{\rm Consider the following simple monic polynomials over $\H:$
   \[
    p_l(z)=z^3-{\bf{k}} z^2+ ({\bf{k-j}})z+({\bf{i+ j}})\,\,\mbox{and}\,\,    p_r(z)=z^3 -z^2{\bf{k}}+ z({\bf{k-j}})+ ({\bf{i+j}}).
    \]
The corresponding companion matrices of  $p_l(z)$ and $p_r(z)$ are given by
$$C_{p_l}= \kbordermatrix{ &1  &  & 2 \\
 2& 0 & \vrule & I \\
  \cline{2-4}
1 & C_{p_l}(3, 1) & \vrule   &  C_{p_l}(3, 2:3)
 }\,\,\mbox{and}\,\,C_{p_r}= C_{p_l}^T,$$
 respectively, 
 where $C_{p_l}(3, 1)= {\bf{-i-j}}$ and $  C_{p_l}(3, 2: 3):=[{\bf{j-k}}, {\bf{k}}].$ Then
\[
C^2_{p_l}=\bmatrix{0 & 0 & 1 \\ \bf{-i-j} & \bf{j-k} & \bf{k} \\ \bf{i-j} & \bf{1-2i-j} & \bf{j-k-1}}\,\,\mbox{and}\,\,
C^2_{p_r}= \bmatrix{0 & \bf{-i-j} & \bf{j-i} \\ 0 & \bf{j-k} & \bf{1-j} \\ 1 & \bf{k} & \bf{j-k-1} }.
\]
This shows that $C_{p_r}^2 \not= (C_{p_l}^2)^T.$
}
\end{exam}
Hence, we can derive results analogous to Theorem \ref{i} for the case of $C_{p_r}^t, t\ge 2.$
\begin{theorem}\label{ii} Consider
$C_{p_r}=\kbordermatrix{ &m-1& & 1\\
 1& 0 &\vrule & C_{p_r}(1, m) \\
  \cline{2-4}
 m-1& I & \vrule & C_{p_r}(2: m, m)
 }$. \\
 {\bf (a)} If $t < m$ is a positive integer, then 
\begin{equation}\label{q1}
 C_{p_r}^t=
 \kbordermatrix{ &m-t& & t\\
 t&  0 &\vrule & C \\
  \cline{2-4}
 m-t& I  & \vrule & D
 },
 \end{equation}
 {\bf (b)} if $t \ge m,$ then
 \begin{eqnarray*}\label{q2}
  C^t_{p_r}=\left[
 \begin{array}{ccccc}
 C^{t-(m-1)}_{p_r}(1: m, m) & C^{t-(m-2)}_{p_r}(1: m, m) &\dots & C^{t-1}_{p_r}(1:m, m) & C_{p_r}^t(1:m, m)
 \end{array}\right]_{m\times m},
 \end{eqnarray*}
 where
 \begin{eqnarray*}
 C &:=& \bmatrix{C_{p_r}(1: t, m) & C^2_{p_r}(1: t, m) & \ldots & C^t_{p_r}(1: t, m)},\\
 D &:=& \bmatrix{C_{p_r}(t+1: m, m) & C^2_{p_r}(t+1: m, m) &\ldots& C^t_{p_r}(t+ 1: m, m)},\\
 C_{p_r}^t(1, m) &:=& C_{p_r}(1, m) \,\,  C_{p_r}^{t-1}(m, m),\,\, \mbox{and}\\
 \,\,C_{p_r}^t(2: m, m) &:=& C_{p_r}^{t-1}(1: m-1, m)+
C_{p_r}(2: m, m) \,\,  C_{p_r}^{t-1}(m, m).
\end{eqnarray*}
\end{theorem}
\proof The proof follows from the proof method of Theorem \ref{i}.$\,\,\,\blacksquare$
%
%
%

 Polynomials from Example \ref{ex5.7} satisfy
   \[
     \tilde{p}_l(z):= \overline{p_r(\overline{z})}= z^3+ {\bf{k}} z^2+({\bf{j-k}})z+ ({\bf{-i-j}}),\,\,\mbox{and}\,\,
     \tilde{p_r}(z):= \overline{p_l(\overline{z})}= z^3+ z^2{\bf{k}}+z({\bf{j-k}})- ({\bf{i+j}}).
     \]
 Thus the companion matrices corresponding to  $\tilde{p}_l(z)$ and $ \tilde{p}_r(z)$ are given by
\[
C_{\tilde{p}_l}= \overline{C_{p_l}} \,\mbox{and}\, C_{\tilde{p}_r}=\overline{C_{p_r}},
\]
respectively.
Next,
\[
C^2_{\tilde{p}_l}= \bmatrix{0 & 0& 1 \\ \bf{i+j} & \bf{-j+k} & \bf{-k} \\ \bf{i-j} & \bf{1+j} & \bf{k-j-1}}\, \mbox{and}\,
C^2_{\tilde p_r}= \bmatrix{0 & \bf{i+ j} & \bf{j- i} \\ 0 & \bf{-j+k} & \bf{1+ 2i+ j} \\ 1 & \bf{-k} & \bf{-1 -j+ k} }.\]
Then
$$ \dm{\max_{ 1 \leq i \leq 3}}\left[ (r_i' (C^2_{p_l}) )^{1/2} \right]=2.3655\,\,\,\mbox{and}\,\,\,
\dm{\max_{ 1 \leq i \leq 3}}\left[ (r_i' (C^2_{\tilde{p_r}} ) )^{1/2} \right]=1.9656,$$
%
%
$$
 \dm{\max_{ 1 \leq i \leq 3}}\left[ \left(r_i'\left(C^2_{p_r}\right) \right)^{1/2}
 \right]=1.9319 \,\,\mbox{and}\,\, \dm{\max_{ 1 \leq i \leq 3}}\left[ (r_i'(C^2_{\tilde{p_l}}) )^{1/2}  \right]=2.1355.
 $$
  Now, we have
\begin{eqnarray*}
  \dm{\max_{ 1 \leq i \leq 3}}\left[ (r_i' (C^2_{p_l}) )^{1/2} \right] &\not =&
  \dm{\max_{ 1 \leq i \leq 3}}\left[ (r_i' (C^2_{\tilde{p_r}} ) )^{1/2} \right]\,\,\mbox{and}\\
     \dm{\max_{ 1 \leq i \leq 3}}\left[ \left(r_i'\left(C^2_{p_r}\right) \right)^{1/2}
 \right] &\not=& \dm{\max_{ 1 \leq i \leq 3}}\left[ (r_i'(C^2_{\tilde{p_l}}) )^{1/2}  \right].
 \end{eqnarray*}

Further, we have the following bounds for the zeros of $p_l(z)$ and $p_r(z)$ for $\gamma \in [0,\, 1].$
\begin{theorem}\label{T4}
 Let $p_l(z)$ and $p_r(z)$ be the simple monic polynomials over $\H$ of degree $m$ and let $C_{p_l}^t$ and
 $C_{p_r}^t\, (t \geq 2$) be the $t$-th power of the
companion matrices $C_{p_l}$ and $C_{p_r},$ corresponding to $p_l(z)$ and $p_{r}(z),$ respectively. Then, for $\gamma \in [0, 1]$
bounds for every zero $\tilde{z}$ of $p_l(z)$ satisfy the following inequalities:
\begin{eqnarray}
\left(\dm{\max_{ 1 \leq i \leq m}}\left[ \left(r_i'\left(C^t_{q_l}\right) \right)^{\gamma/t}
 \left(c_i'\left(C^t_{q_l}\right) \right)^{(1-\gamma)/t}\right]\right)^{-1}
 \le
 |\tilde{z}| \le \dm{\max_{ 1 \leq i \leq m}}\left[ \left(r_i'\left(C^t_{p_l}\right) \right)^{\gamma/t}
 \left(c_i'\left(C^t_{p_l}\right) \right)^{(1-\gamma)/t}\right],\label{neqn3}
 \end{eqnarray}
\begin{eqnarray}
\left(\dm{\max_{ 1 \leq i \leq m}}\left[ \left(r_i'\left(C^t_{\tilde{q_r}}\right) \right)^{\gamma/t}
 \left(c_i'\left(C^t_{\tilde{q_r}}\right) \right)^{(1-\gamma)/t}\right]\right)^{-1}
 \le
 |\tilde{z}| \le \dm{\max_{ 1 \leq i \leq m}}\left[ \left(r_i'\left(C^t_{\tilde{p_r}}\right) \right)^{\gamma/t}
 \left(c_i'\left(C^t_{\tilde{p_r}}\right) \right)^{(1-\gamma)/t}\right], \label{neqn4}
 \end{eqnarray}
%
 and bounds for every zero $\tilde{z}$ of $p_r(z)$ satisfy the following inequalities:
\begin{eqnarray}
\left(\dm{\max_{ 1 \leq i \leq m}}\left[ \left(r_i'\left(C^t_{q_r}\right) \right)^{\gamma/t}
 \left(c_i'\left(C^t_{q_r}\right) \right)^{(1- \gamma)/t}\right]\right)^{-1}
 \le
 |\tilde{z}|
 \le \dm{\max_{ 1 \leq i \leq m}}\left[ \left(r_i'\left(C^t_{p_r}\right) \right)^{\gamma/t}
 \left(c_i'\left(C^t_{p_r}\right) \right)^{(1- \gamma)/t}\right],\label{neqn5}
\end{eqnarray}
\begin{eqnarray}
\left(\dm{\max_{ 1 \leq i \leq m}}\left[ \left(r_i'\left(C^t_{\tilde{q_l}}\right) \right)^{\gamma/t}
 \left(c_i'\left(C^t_{\tilde{q_l}}\right) \right)^{(1-\gamma)/t}\right]\right)^{-1}
 \le
 |\tilde{z}| \le \dm{\max_{ 1 \leq i \leq m}}\left[ \left(r_i'\left(C^t_{\tilde{p_l}}\right) \right)^{\gamma/t}
 \left(c_i'\left(C^t_{\tilde{p_l}}\right) \right)^{(1- \gamma)/t}\right] \label{neqn6}.
 \end{eqnarray}
\end{theorem}
\proof Let $\lam$ be a left eigenvalue of $C_{p_l}.$
Then by Proposition \ref{p2}, $\lam^t$ ( $t \geq 2$ is positive integer) is a left eigenvalue
of $C^t_{p_l}.$ Hence by applying Theorem \ref{os1}, we get (\ref{neqn3}).

\noin By Lemma \ref{prop4}, $\overline{\lam}$ is a left eigenvalue of $C_{\tilde{p_r}}$ and by
Proposition \ref{p3}, $(\overline{\lam})^t$ is a left eigenvalue of $(C_{\tilde{p_r}})^t.$ Then from Theorem \ref{os1}, (\ref{neqn4}) follows.

\noin The proof of (\ref{neqn5}) and (\ref{neqn6}) are similar$.\,\,\,\blacksquare$

Substituting $t=2$ and $\gamma=1$ in Theorem \ref{T4}, we have the following corollary.
\begin{corollary}\label{pc}
 Let $p_l(z)$ and $p_r(z)$ be the simple monic polynomials over $\H$ of degree $m.$ Then
 bounds for every zero $\tilde{z}$ of $p_l(z)$ satisfy the following inequalities:
\begin{eqnarray}
        \frac{1}{\beta_1}    \leq |\tilde{z}| \leq \alpha_1, 
        \end{eqnarray}
\begin{eqnarray}
\frac{1}{\beta_2}    \leq |\tilde{z}| \leq \alpha_2, 
\end{eqnarray}
where
\begin{eqnarray*}
\alpha_1 &=& \dm{\max \left\{ 1, \left( \sum_{j=0}^{m-1} |q_j|   \right)^{1/2},
\left( \sum_{j=0}^{m-1}|q_{m-1} q_j - q_{j-1}|   \right)^{1/2} \right\}},\\
\alpha_2 &=& \dm{ \max_{ 2 \leq j \leq m-1 } \left\{ \left( |q_0| + |\overline{q_0}\,\, \overline{q_{m-1}}|   \right)^{1/2},
\left(   |q_1| + |\overline{q_1}\,\, \overline{q_{m-1}} - \overline{q_0}|    \right)^{1/2},
\left(  1+ |q_j| + |\overline{q_j}\,\, \overline{q_{m-1}} - \overline{q_{j-1}}|    \right)^{1/2}   \right\}},\\
\beta_1 &=&  \dm{\max \left\{ 1, \left( \sum_{j=1}^{m-1} |q^{-1}_0 q_j|   \right)^{1/2},
\left( \sum_{j=0}^{m-1}|q^{-1}_0 q_1 q^{-1}_0q_{m-j} -  q^{-1}_0 q_{m-j+1}|   \right)^{1/2} \right\}},\\
\beta_2 &=&  \max_{ 2 \leq j \leq m-1 }  \bigg\{  \left(  |q^{-1}_0| + |\overline{q^{-1}_0}\,\,\, \overline{q_1 q^{-1}_0}|
\right)^{1/2},
\left(  |q_{m-1} q^{-1}_0| + |\overline{ q_{m-1} q^{-1}_0}\,\,\, \overline{q_1 q^{-1}_0}- \overline{q^{-1}_0}| \right)^{1/2},\\ &&
\left(  1+ |q_{m-j} q^{-1}_0| + |\overline{ q_{m-j} q^{-1}_0}\,\,\, \overline{q_1 q^{-1}_0}- \overline{q_{m-j+ 1}q^{-1}_0}|
\right)^{1/2} \bigg\},
\end{eqnarray*}
%
%
 and bounds for every zero $\tilde{z}$ of $p_r(z)$ satisfy the following inequalities:
\begin{eqnarray}      
\frac{1}{\beta_3}    \leq |\tilde{z}| \leq \alpha_3, 
 \end{eqnarray}
\begin{eqnarray}
\frac{1}{\beta_4}    \leq |\tilde{z}| \leq \alpha_4, 
\end{eqnarray}
where
\begin{eqnarray*}
\alpha_3 &=& \dm{ \max_{ 2 \leq j \leq m-1 } \left\{ \left( |q_0| + |q_0\,\, q_{m-1}|   \right)^{1/2},
\left(   |q_1| + |q_1\,\, q_{m-1} - q_0|    \right)^{1/2} ,
\left(  1+ |q_j| + |q_j\,\, q_{m-1} - q_{j-1}|    \right)^{1/2}   \right\}},\\
\alpha_4 &=& \dm{\max \left\{ 1, \left( \sum_{j=0}^{m-1} |q_j|   \right)^{1/2},
\left( \sum_{j=0}^{m-1}|\overline{q_{m-1}} \,\, \overline{q_j} - \overline{q_{j-1}}|   \right)^{1/2} \right\}},\\
\beta_3 &=&  \max_{ 2 \leq j \leq m-1 }  \bigg\{ \small{ \left(  |q^{-1}_0| + |q^{-1}_0\,\,\, q_1 q^{-1}_0|
\right)^{1/2}},
\left(  |q_{m-1} q^{-1}_0| + |q_{m-1} q^{-1}_0\, q_1 q^{-1}_0- q^{-1}_0| \right)^{1/2},\\ &&
\left(  1+|q_{m-j} q^{-1}_0| + |q_{m-j} q^{-1}_0\, q_1 q^{-1}_0- q_{m-j+1}q^{-1}_0|
\right)^{1/2}\bigg\},\\
\beta_4 &=& \dm{\max \left\{ 1, \left( \sum_{j=1}^{m-1} |q^{-1}_0 q_j|   \right)^{1/2},
\left( \sum_{j=0}^{m-1}|\overline{q^{-1}_0 q_1} \,\,\,\, \overline{q^{-1}_0q_{m-j}} -  \overline{q^{-1}_0 q_{m-j+1}}|   \right)^{1/2} \right\}},
q_{-1}=0=q_{m+1}, q_m=1.
\end{eqnarray*}
\end{corollary}
\proof The proof follows from Theorem \ref{T4} and Appendix \ref{AP1}.$\,\,\,\blacksquare$
\begin{exam}\label{e1}{\rm
  Consider the following polynomials $p_l(z)$ and $p_r(z)$ over $\H$:
  \[p_l(z)=z^6+ ({\bf{i}}+ 3{\bf{k}}) z^5+ (3+ {\bf{j}})z^4+(5{\bf{i}}+ 15{\bf{k}}) z^3+ (-4+ 5{\bf{j}})z^2+ (6{\bf{i}}+ 18{\bf{k}})z+ (
  6{\bf{j}}  -12),\]
\[p_r(z)=z^6+ z^5 ({\bf{i}}+ 3{\bf{k}})+ z^4 (3+ {\bf{j}})+z^3 (5{\bf{i}}+ 15{\bf{k}})+ z^2 (-4+5{\bf{j}})+z (6{\bf{i}}+18{\bf{k}})+ (
  6{\bf{j}}  -12).\]}
 \end{exam}
 The zeros of $p_l(z)$ are given in \cite{rp01}. Moreover, we find the zeros of $p_r(z)$ by Niven's algorithm \cite{in41}.

\begin{table}[!h]
    \caption{Zeros and bounds for the zeros of $p_l(z)$ and $p_r(z).$}
    \begin{subtable}{.6\linewidth}
      \centering
       \caption{Zeros of $p_l(z)$ and $p_r(z)$  and their absolute values.}\label{t}
               \begin{tabular}{|cccc|}\hline
$ z_1$ & $ |z_1|$ & $z_2$ & $ |z_2|$\\
\hline
   $-{\bf{i}}-2{\bf{k}}$ &  $  2. 2361$ &$-0.4{\bf{i}}-2.2{\bf{k}}$ & $  2. 2361$ \\
  $[{\bf{i}}\sqrt{3}]$ & $1. 7321$ & $[{\bf{i}}\sqrt{3}]$& $1. 7321$\\ 
  $ [{\bf{i}}\sqrt{2}]$ & $ 1. 4142$ & $ [{\bf{i}}\sqrt{2}]$ & $ 1. 4142$ \\
  $-0.6{\bf{i}}-0. 8{\bf{k}}$ & $ 1$  &$-{\bf{k}}$ & $ 1$ \\
  \hline 
\end{tabular}
    \end{subtable}%
    \begin{subtable}{.4\linewidth}
      \centering
        \caption{Lower and upper bounds for the zeros of $p_l(z)$ and $p_r(z).$}\label{t4}
        \begin{tabular}{|ccc|}\hline
Example \ref{e1}   & lower bound  &  upper bound \\ 
   \hline
Corollary \ref{co} (1) & $0. 4142 $      & $19. 9737$   \\ 
Corollary   \ref{co} (2) &$0. 2766 $      & $60. 9291$   \\ 
  Theorem \ref{u1}, $\gamma = 1/4$ &$0. 3744 $      & $8. 1415$   \\ 
   \hline   
\end{tabular}
    \end{subtable} 
\end{table}
where
$z_1:= $ the set of zeros of $p_l(z),$ 
$z_2:= $ the set of zeros of $p_r(z)$
\begin{table}[h!]
 \begin{center}
\begin{tabular}{ccccc}\hline
Example \ref{ex5.7}  &&  lower bound &&   lower bound\\ 

   \hline
Corollary \ref{pc} $1(a)$ && $ 0. 6156 $  && $2. 3655$   \\ 
   Corollary \ref{pc} $1(b)$  && $ 0. 6078 $ && $1. 9656$   \\ 
   Corollary \ref{pc} $2(a)$  && $0. 6078$ && $ 1. 9319 $   \\ 
  Corollary \ref{pc} $2(b)$  && $0. 6436$  && $2. 1355$   \\ 
   \hline
\end{tabular}
 \caption{Lower and upper bounds for the zeros of $p_l(z)$ and $p_{r}(z).$}
\end{center}
\end{table}

\section{Conclusion}\label{cf}
In this paper, we have derived Ostrowski type theorem for left eigenvalues of a quaternionic matrix that generalizes
Ostrowski type theorem for right eigenvalues of a quaternionic matrix when all the diagonal entries of a quaternionic matrix
are real. We have derived a corrected version of the Brauer type theorem for left eigenvalues for the
deleted absolute column sums of a quaternionic matrix. Moreover, we have extended
localization theorems by applying the generalized H$\ddot{\mbox{o}}$lder inequality for left as well as right eigenvalues of 
a quaternionic matrix. 
Bounds for the zeros of quaternionic polynomials have derived.
As a consequence, we have shown that some of our bounds are sharper than the bound given in  \cite{go09}. Further, 
we have derived bounds via the powers of companion matrices which are always sharper than the bound given in \cite{go09}.
\vone

\noin{\bf Acknowledgements:}
The authors would like to thank the reviewer and editor for their valuable comments and suggestions to improve the manuscript. They also thank Professor  Ivan Slapni$\check{\mbox{c}}$ar for careful reading and helpful comments for the improvement of the manuscript.

\begin{appendix} \section{Appendix}\label{AP1}
In this appendix, we state formulas for the squares of quaternionic companion matrices. For $t= 2$,  Theorem \ref{i} implies
$$
 C_{p_l}^2=
 \kbordermatrix{ &2& & m-2\\
 m-2&  0 &\vrule & I \\
  \cline{2-4}
 2& C  & \vrule & D
 }, \,\mbox{where}\,\,
C:=\bmatrix{C_{p_l}(m,1: 2)\\C^2_{p_l}(m,1: 2) }=\bmatrix{-q_0 & -q_1 \\ q_{m-1} q_0 & q_{m-1} q_1-q_0}\,\,$$
 and
 $$D= \bmatrix{C_{p_l}(m,3: m)\\C^2_{p_l}(m, 3: m)}=\bmatrix{-q_2 & -q_3 & \ldots & -q_{m-1} \\
 q_{m-1} q_2-q_1 &  q_{m-1} q_3-q_1 &\ldots & (q_{m-1})^2-q_{m-2}},
$$
$$
C_{ \tilde{p_l}}^2=
 \kbordermatrix{ &2& & m-2\\
 m-2&  0 &\vrule & I \\
  \cline{2-4}
 2& C  & \vrule & D
 }, \,\mbox{where}\,\,
C= \bmatrix{C_{ \tilde{p_l}}(m,1:2)\\C^2_{ \tilde{p_l}}(m,1:2) }= \bmatrix{-\overline{q_0} & -\overline{q_1} \\ \overline{q_{m-1}}\,\,
\overline{q_0} & \overline{q_{m-1}} \, \,\overline{q_1}-\overline{q_0}}$$
and
 $$D=\bmatrix{C_{ \tilde{p_l}}(m,3:m)\\C^2_{ \tilde{p_l}}(m, 3:m)}=\bmatrix{-\overline{q_2} & -\overline{q_3}
 & \ldots & -\overline{q_{m-1}} \\
 \overline{q_{m-1}} \,\,  \overline{q_2}-\overline{q_1} &  \overline{q_{m-1}} \,\, \overline{q_3}-\overline{q_1}
 &\ldots & (\overline{q_{m-1}})^2-\overline{q_{m-2}}},
$$
$$
 C_{q_l}^2=
 \kbordermatrix{ &2& & m-2\\
 m-2&  0 &\vrule & I \\
  \cline{2-4}
 2& C  & \vrule & D
 }, \,\mbox{where}\,\,
C=\bmatrix{-q^{-1}_0 & -q^{-1}_0 q_{m-1} \\  q^{-1}_0 q_1 q^{-1}_0 & q^{-1}_0 q_1 q^{-1}_0 q_{m-1} -q^{-1}_0}\,\,$$
and
 $$D=\bmatrix{- q^{-1}_0 q_{m-2} & \ldots &  - q^{-1}_0 q_1 \\
 q^{-1}_0 q_1 q^{-1}_0 q_{m-2}- q^{-1}_0 q_{m-1} & \ldots & (q^{-1}_0 q_1)^2-q^{-1}_0q_2    },
$$
$$
 C_{\tilde{q_l}}^2=
 \kbordermatrix{ &2& & m-2\\
 m-2&  0 &\vrule & I \\
  \cline{2-4}
 2& C  & \vrule & D
 }, \,\mbox{where}\,\,
 C=\bmatrix{-\overline{q^{-1}_0} & -\overline{q^{-1}_0 q_{m-1}} \\
&   \\ \overline{q^{-1}_0 q_1}\,\,\overline{ q^{-1}_0} &
\overline{q^{-1}_0 q_1}\,\, \overline{q^{-1}_0 q_{m-1}} -\overline{q^{-1}_0}}\,\,$$
and
 $$D=\bmatrix{- \overline{q^{-1}_0 q_{m-2}} & \ldots &  - \overline{q^{-1}_0 q_1} \\
 & \\
 \overline{q^{-1}_0 q_1}\,\,\,\overline{ q^{-1}_0 q_{m-2}}- \overline{q^{-1}_0 q_{m-1}} & \ldots &
 \left(\overline{q^{-1}_0 q_1}\right)^2-\overline{q^{-1}_0q_2 }   }.
$$

 For $t=2$, Theorem \ref{ii} implies
$$
C_{p_r}^2=
 \kbordermatrix{ &m-2& & 2\\
 2&  0 &\vrule & C \\
  \cline{2-4}
 m-2& I  & \vrule & D
 }, \mbox{where}\,\,
C=\bmatrix{C_{p_r}(1:2, m) & C^2_{p_r}(1:2, m) }=\bmatrix{-q_0 & q_0 q_{m-1} \\ -q_1 & q_1 q_{m-1}-q_0   },$$ and
 $$D=\bmatrix{C_{p_r}(3:m, m) & C^2_{p_r}(3:m, m)}= \bmatrix{ -q_2 & q_2 q_{m-1}-q_1 \\ -q_3 & q_3 q_{m-1}-q_2 \\
 \vdots & \vdots \\
 -q_{m-1}  & (q_{m-1})^2 - q_{m-2}},
 $$
$$
C_{\tilde{p_r}}^2=
 \kbordermatrix{ &m-2& & 2\\
 2&  0 &\vrule & C \\
  \cline{2-4}
 m-2& I  & \vrule & D
 }, \mbox{where}\,\,
C=\bmatrix{-\overline{q_0} & \overline{q_0}\,\, \overline{q_{m-1}} \\
-\overline{q_1} & \overline{q_1}\,\, \overline{q_{m-1}}-\overline{q_0}   }\,\, \mbox{and}\,\,
  D=\bmatrix{ -\overline{q_2} & \overline{q_2}\,\, \overline{q_{m-1}}-\overline{q_1} \\ -\overline{q_3} & \overline{q_3}\,\,
 \overline{q_{m-1}}-\overline{q_2} \\
 \vdots & \vdots \\
 -\overline{q_{m-1}}  & \left(\overline{q_{m-1}}\right)^2 - \overline{q_{m-2}}},
 $$
$$
C_{q_r}^2=
 \kbordermatrix{ &m-2& & 2\\
 2&  0 &\vrule & C \\
  \cline{2-4}
 m-2& I  & \vrule & D
 }, \mbox{where}
$$
$$C=\bmatrix{-q^{-1}_0 & q^{-1}_0\,\, q_{1} q^{-1}_0 \\
-q_{m-1} q^{-1}_0 & q_{m-1} q^{-1}_0 q_1 q^{-1}_0-q^{-1}_0   }\,\,\mbox{and}\,\,
 D=\bmatrix{ -q_{m-2} q^{-1}_0 & q_{m-2} q^{-1}_0 q_1 q^{-1}_0-q_{m-1} q^{-1}_0  \\
 \vdots & \vdots \\
 -q_{1} q^{-1}_0 & (q_1 q^{-1}_0)^2 - q_{2} q^{-1}_0},
 $$
$$
C_{\tilde{q_r}}^2=
 \kbordermatrix{ &m-2& & 2\\
 2&  0 &\vrule & C \\
  \cline{2-4}
 m-2& I  & \vrule & D
 }, \mbox{where}
$$
$$C=\bmatrix{-\overline{q^{-1}_0} & \overline{q^{-1}_0}\,\,\, \overline{q_{1} q^{-1}_0} \\
& \\
-\overline{q_{m-1} q^{-1}_0} & \overline{q_{m-1} q^{-1}_0}\,\, \overline{q_1 q^{-1}_0}-\overline{q^{-1}_0}   }\,\,\mbox{and}\,\,
 D=\bmatrix{ -\overline{q_{m-2} q^{-1}_0} & \overline{q_{m-2} q^{-1}_0}\,\,\, \overline{q_1 q^{-1}_0}-\overline{q_{m-1} q^{-1}_0}  \\
 \vdots & \vdots \\
 -\overline{q_{1} q^{-1}_0} & \left(\overline{q_1 q^{-1}_0}\right)^2 - \overline{q_{2} q^{-1}_0}}.
 $$
\end{appendix}


\begin{thebibliography}{10}

\bibitem{sla95}
S.L. Adler.
\newblock {\em Quaternionic {Q}uantum {M}echanics and {Q}uantum {F}ields}.
\newblock Oxford University Press, New York, 1995.

\bibitem{a99}
A.~Baker.
\newblock Right eigenvalues for quaternionic matrices{:} a topological
  approach.
\newblock {\em Linear Algebra Appl.}, 286:303--309, 1999.

\bibitem{a46}
A.~Brauer.
\newblock Limits for the characteristic roots of a matrix.
\newblock {\em Duke Math. J.}, 13:387--395, 1946.

\bibitem{jd02}
J.~H. Conway and D.~A. Smith.
\newblock {\em On Quaternions and Octonions: Their Geometry, Arithmetic, and
  Symmetry}.
\newblock A K Peters Natick, 2002.

\bibitem{s31}
S.~Ger$\check{\mbox{s}}$gorin.
\newblock $\ddot{\mbox{u}}$ber die abgrenzung der eigenwerte einer matrix.
\newblock {\em Izv. Akad. Nauk SSSR Ser. Mat.}, 1:749--754, 1931.

\bibitem{arv89}
A.B. Gerstner, R.~Byers, and V.~Mehrmann.
\newblock A quaternion {QR} algorithm.
\newblock {\em Numerih. Mathek.}, 55:83--95, 1989.

\bibitem{bt65}
B.~Gordon and T.~S. Motzkin.
\newblock On the zeros of polynomials over division rings.
\newblock {\em Trans. Amer. Math. Soc.}, 116:218--226, 1965.

\bibitem{t80}
T.L. Hankins.
\newblock {\em Sir William Rowan Hamilton}.
\newblock The Johns Hopkins University Press, Baltimore, 1980.

\bibitem{m04}
M.~A. Hassan.
\newblock Inequalities and bounds for the zeros of polynomials using
  {P}erron-{F}robenius and {G}erschgorin theories.
\newblock {\em Proceedings of American Control Conference, Boston,
  Massachusetts. Evanston, IL: American Automatic Control Council},
  3:2745--2750, 2004.

\bibitem{nf03}
N.J. Higham and F.~Tisseur.
\newblock Bounds for eigenvalues of matrix polynomials.
\newblock {\em Linear Algebra Appl.}, 358:5--22, 2003.

\bibitem{rc96}
R.~A. Horn and C.~R. Johnson.
\newblock {\em Matrix Analysis}.
\newblock Cambridge University Press, New York, 1996.

\bibitem{rf12}
R.~A. Horn and F.~Zhang.
\newblock A generalization of the complex autonne-takagi factorization to
  quaternion matrices.
\newblock {\em Linear and Multilinear Algebra}, 60:1239--1244, 2012.

\bibitem{lw01}
L.~Huang and W.~So.
\newblock On left eigenvalues of a quaternionic matrix.
\newblock {\em Linear Algebra Appl.}, 323:105--116, 2001.

\bibitem{ddg10}
D.~Janovsk$\acute{\mbox{a}}$ and G.~Opfer.
\newblock The classification and the computation of the zeros of quaternionic,
  two-sided polynomials.
\newblock {\em Numerih. Mathek.}, 115:81--100, 2010.

\bibitem{dg10}
D.~Janovsk$\acute{\mbox{a}}$ and G.~Opfer.
\newblock A note on the computation of all zeros of simple quaternionic
  polynomials.
\newblock {\em {SIAM} J. Numer. Anal.}, 48:244--256, 2010.

\bibitem{wzcl08}
W.~Junliang, Z.~Limin, C.~Xiangping, and L.~Shengjie.
\newblock The estimation of eigenvalues of sum, difference, and tensor product
  of matrices over quaternion division algebra.
\newblock {\em Linear Algebra Appl.}, 428:3023--3033, 2008.

\bibitem{gp02}
G.~Kamberov, P.~Norman, F.~Pedit, and U.~Pinkall.
\newblock {\em Quaternions, Spinors, and Surfaces, Contemporary Mathematics,
  vol. 299}.
\newblock Amer. Math. Soc., Province, 2002.

\bibitem{m99}
M.~Karow.
\newblock Self-adjoint operators and pairs of {H}ermitian forms over the
  quaternions.
\newblock {\em Linear Algebra Appl.}, 299:101--117, 1999.

\bibitem{l49}
H.C. Lee.
\newblock Eigenvalues of cannonical forms of matrices with quaternion
  coefficients.
\newblock {\em Pro. R. I. A. 52}, Sec A:253--260, 1949.

\bibitem{sgv06}
S.~D. Leo, G.~Ducati, and V.~Leonardi.
\newblock Zeros of unilateral quaternionic polnomials.
\newblock {\em El. J. Lin. Alg.}, 15:297--313, 2006.

\bibitem{in41}
I.~Niven.
\newblock Equations in quaternions.
\newblock {\em Amer. Math. Monthly}, 48:654--661, 1941.

\bibitem{go09}
G.~Opfer.
\newblock Polynomials and vandermonde matrices over the field of quaternions.
\newblock {\em Electr. Trans. Num. Anal.}, 36:9--16, 2009.

\bibitem{am37}
A.~M. Ostrowski.
\newblock $\ddot{U}$ber die determinanten mit uberwiegender hauptdiagonale.
\newblock {\em Comment. Math. Helv.}, 10:69--96, 1937b.

\bibitem{rpp08}
R.~Pereira and P.~Rocha.
\newblock On the determinant of quaternionic polynomial matrices and its
  application to system stability.
\newblock {\em Math. Meth. Appl. Sci.}, 31:99--122, 2008.

\bibitem{rppv05}
R.~Pereira, P.~Rocha, and P.~Vettori.
\newblock Algebraic tools for the study of quaternionic behavioral systems.
\newblock {\em Linear Algebra Appl.}, 400:121--140, 2005.

\bibitem{am04}
A.~Pogorui and M.~Shapiro.
\newblock On the structure of the set of zeros of quaternionic polynomials.
\newblock {\em Complex Var. and Elliptic Funct.}, 49:379--389, 2004.

\bibitem{l08}
L.~Rodman.
\newblock Pairs of hermitian and skew hermitian quaternionic matrices canonical
  forms and their applications.
\newblock {\em Linear Algebra Appl.}, 429:981--1019, 2008.

\bibitem{lr12}
L.~Rodman.
\newblock Stability of invariant subspaces of quaternion matrices.
\newblock {\em Complex. Anal. Oper.Theory}, 6:1069--1119, 2012.

\bibitem{l14}
L.~Rodman.
\newblock {\em Topics in {Q}uaternion {L}inear {A}lgebra}.
\newblock Princeton University Press, Princeton, {NJ}, 2014.

\bibitem{rp01}
R.~Ser$\hat{\mbox{o}}$dio, E.~Pereira, and J.~Vit$\acute{\mbox{o}}$ria.
\newblock Computing the zeros of quaternion polynomials.
\newblock {\em Comp. Math. Appl.}, 42:1229--1237, 2001.

\bibitem{cd11}
C.~C. Took and D.~P. Mandic.
\newblock Augmented second-order statistics of quaternion random signals.
\newblock {\em Signal Processing}, 91:214--224, 2011.

\bibitem{cdf11}
C.~C. Took, D.~P. Mandic, and F.~Zhang.
\newblock On the unitary diagonalisation of a special class of quaternion
  matrices.
\newblock {\em Applied Mathematics Letters}, 24:1806--1809, 2011.

\bibitem{r04}
R.S. Varga.
\newblock {\em Ger$\check{\mbox{s}}$gorin and {H}is {C}ircles}.
\newblock Springer-Verlag, Berlin, 2004.

\bibitem{jw08}
J.~Wu.
\newblock Distribution and estimation for eigenvalues of real quaternion
  matrices.
\newblock {\em Comp. Math. Appl.}, 55:1998--2004, 2008.

\bibitem{fz97}
F.~Zhang.
\newblock Quaternions and matrices of quaternions.
\newblock {\em Linear Algebra Appl.}, 251:21--57, 1997.

\bibitem{fz07}
F.~Zhang.
\newblock Ger$\check{\mbox{s}}$gorin type theorems for quaternionic matrices.
\newblock {\em Linear Algebra Appl.}, 424:139--155, 2007.

\bibitem{lyj12}
L.~Zou, Y.~Jiang, and J.~Wu.
\newblock Location for the right eigenvalues of quaternion matrices.
\newblock {\em J. Appl. Math. Comput.}, 38:71--83, 2012.

\end{thebibliography}
\end{document}